\documentclass[12pt]{amsart} 
\usepackage{amssymb}
\usepackage[colorlinks=true]{hyperref}
\usepackage{graphicx}

\newcommand{\bbK}{\mathbb{K}}

\newcommand{\bbQ}{\mathbb{Q}}

\newcommand{\bbZ}{\mathbb{Z}}

\newcommand{\cN}{\mathcal{N}}
\newcommand{\cO}{\mathcal{O}}

\newcommand{\fa}{\mathfrak{a}}
\newcommand{\fp}{\mathfrak{p}}
\newcommand{\fq}{\mathfrak{q}}

\DeclareMathOperator{\core}{sf}

\numberwithin{equation}{section}

\newtheorem{theorem}{Theorem}[section]
\newtheorem{conjecture}[theorem]{Conjecture}
\newtheorem{corollary}[theorem]{Corollary}
\newtheorem{lemma}[theorem]{Lemma}
\newtheorem{proposition}[theorem]{Proposition}

\theoremstyle{definition}
\newtheorem{example}{Example}[section]
\newtheorem*{remark-nonum}{Remark}
\newtheorem*{remarks-nonum}{Remarks}

\newcommand{\bSqrThm}{1.4} 

\newcommand{\dioSect}{2}
\newcommand{\lemDio}{2.1}
\newcommand{\lemHypg}{2.2}

\newcommand{\eqNdn}{2.1}       
\newcommand{\eqGDef}{2.4}      
\newcommand{\eqApproxDef}{2.5} 
\newcommand{\eqDPrimeDef}{2.6} 
\newcommand{\eqKZero}{2.8}     

\newcommand{\repProp}{3.1}
\newcommand{\lbLemma}{3.5}
\newcommand{\gapLemma}{3.8}
\newcommand{\miscLemma}{3.9}
\newcommand{\gnLemma}{3.12}   
\newcommand{\EQLBLemma}{3.13} 

\newcommand{\eqZetaIJRel}{3.23} 
\newcommand{\gpNumEqn}{3.24} 
\newcommand{\gpEqnb}{3.28}   

\newcommand{\genPropSection}{4}
\newcommand{\genProp}{4.1}
\newcommand{\subSectPreq}{4.1} 
\newcommand{\rEqualSubsection}{4.3} 
\newcommand{\subsectStepiii}{4.4}

\newcommand{\yEllLBEqnA}{4.2} 
\newcommand{\eqELB}{4.5}      
\newcommand{\eqQLB}{4.6}      
\newcommand{\eqQUB}{4.7}      
\newcommand{\eqEllUB}{4.8}    
\newcommand{\eqQAbs}{4.9}     
\newcommand{\eqYEllUBa}{4.11} 
\newcommand{\eqYEllUB}{4.15}  
\newcommand{\eqnYKUBa}{4.19}  
\newcommand{\eqnYKUBb}{4.22}  

\begin{document}

\title[Sharp bounds on \ldots]{Sharp bounds on the number of squares in recurrence sequences and solutions of $X^{2}-\left( a^{2}+b \right) Y^{4}=-b$}

\author{Paul M Voutier}
\address{London, UK}
\email{Paul.Voutier@gmail.com}

\date{}

\begin{abstract}
We obtain best possible results for the number of coprime positive integer
solutions of the equation in the title
when $a$ is a positive integer, $b=p^{m}$, $2p^{m}$ or $4p^{m}$, where
$m$ is a non-negative integer, $p$ is prime, $\gcd \left( a^{2}, b \right)$ is
squarefree and $X^{2}- \left( a^{2}+b \right) Y^{2}=-4$
has a solution in positive integers.
We prove our results by establishing best possible bounds for the number of
distinct squares in certain binary recurrence sequences, including those
associated with such equations.
\end{abstract}


\maketitle

\section{Introduction}

\subsection{Background}

The study of integer solutions of the diophantine equations $aX^{2}-bY^{4}=c$
where $a, b$ are positive integers and $c$ is a non-zero integer goes back
over 100 years. Results of Landau and Ostrowski \cite{LO} and, slightly later,
Mordell \cite{Mor1} imply that there are only finitely many integer
solutions of these equations.

For $c=\pm 1, \pm 2, \pm 4$, we have best possible bounds for the number
of positive integer solutions of such equations \cite{Akh,L3,LY}.
Such results started with the work
of Ljunggren and since his work, these equations have been the
subject of much attention (see \cite{Akh} and the references there).

For other values of $c$, only a few results are known (see, for example, \cite{SWY}
mentioned below). The author is unaware of any best possible results for such
values of $c$ prior to this work. The study of these equations for such values
of $c$ appears to be much more difficult.

In \cite{V5}, we presented a technique for bounding the number of distinct
integer squares in certain binary recurrence sequences. If there is only one
family of solutions of the Pell-like equation, $aX^{2}-bY^{2}=c$, and our technique
applies, then we can also bound the number of solutions of $aX^{2}-bY^{4}=c$.

In 2009, Stoll, Walsh and Yuan \cite{SWY} showed that for any non-negative
integer $m$, there are at most three solutions in odd positive integers to
\[
X^{2} - \left( 1+2^{2m} \right) Y^{4} = -2^{2m}.
\]

Here we apply our technique to both generalise and sharpen this result of Stoll,
Walsh and Yuan. In fact, under the much more general conditions considered in
Corollary~\ref{cor:1.2-eqn1} below, we actually obtain the best
possible bounds for the number of coprime solutions in positive integers of
\begin{equation}
\label{eq:2}
X^{2} - \left( a^{2}+b \right) Y^{4} = -b.
\end{equation}

\subsection{Results}
\label{subsect:results}

In the remainder of this paper, we let $a$ and $b$ be positive integers such
that $a^{2}+b$ is not a square and let $\varepsilon=\left( t+u\sqrt{a^{2}+b} \right)/2$
be a unit in $\cO_{\bbQ \left( \sqrt{a^{2}+b} \right)}$ with $t$ and $u$ positive
integers. We will denote by $N_{\varepsilon}$, the norm of $\varepsilon$, $\left( t^{2}-\left( a^{2}+b \right)u^{2} \right)/4$.

We define the two sequences $\left( x_{k} \right)_{k=-\infty}^{\infty}$
and $\left( y_{k} \right)_{k=-\infty}^{\infty}$ by
\begin{equation}
\label{eq:yk-defn}
x_{k}+y_{k} \sqrt{a^{2}+b}
=\left( a + \sqrt{a^{2}+b} \right) \varepsilon^{2k}.
\end{equation}

Observe that $x_{0}=a$, $y_{0}=1$, $y_{1}= \left( (t+au)^{2}+bu^{2} \right)/4$,
$y_{-1}= \left( (t-au)^{2}+bu^{2} \right)/4$ and that both sequences satisfy
the recurrence relation
\begin{equation}
\label{eq:yk-recurrence}
u_{k+1}= \frac{t^{2}+\left( a^{2}+b \right) u^{2}}{2} u_{k}-u_{k-1},
\end{equation}
for all $k \in \bbZ$.

In the notation of \cite{V5}, $d$ there is $a^{2}+b$ here, $b$ there is $1$ here
and so $N_{\alpha}$ there is $-b$ here.

In \cite{V5}, we presented the following conjecture, as well as some supporting
evidence for it.

\begin{conjecture}
\label{conj:1-seq}
Let $a$ and $b$ be two positive integers such that $a^{2}+b$
is not a square and let $\left( y_{k} \right)_{k=-\infty}^{\infty}$
be defined by \eqref{eq:yk-defn}.

\noindent
{\rm (a)} There are at most three distinct integer squares among the $y_{k}$'s.

\noindent
{\rm (b)} Moreover, if $b$ is a square, then there are at most two
distinct integer squares among the $y_{k}$'s
\end{conjecture}

We have been able to prove this conjecture under certain conditions on $a$ and
$b$. See Theorem~\bSqrThm{} of \cite{V5}, where we prove part~(b) of the above
conjecture for all but a very specific set of sequences. Here we prove this
conjecture when $b$ is a prime power, or nearly a prime power.

\begin{theorem}
\label{thm:1.2-seq}
Let $a$, $m$ and $p$ be non-negative integers with $a \geq 1$ and $p$ a prime.
Conjecture~$\ref{conj:1-seq}$ holds when $b=p^{m}$, $2p^{m}$, $4p^{m}$, $8p^{m}$
or $16p^{m}$.
\end{theorem}

\begin{corollary}
\label{cor:1.2-eqn1}
Let $a$, $m$ and $p$ be non-negative integers with $a \geq 1$ and $p$ a prime
and put $b=p^{m}$, $2p^{m}$ or $4p^{m}$. Suppose
that $\gcd \left( a^{2}, b \right)$ is squarefree, $a^{2}+b$ is not a square,
$b=p^{m}$ if $a^{2}+b \equiv 1 \pmod{8}$, and
that $X^{2}-\left( a^{2}+b \right)Y^{2}=-4$ has a solution in positive integers.

\noindent
{\rm (a)} If $b$ is a square, then \eqref{eq:2} has at most two coprime positive
integer solutions.

\noindent
{\rm (b)} If $b$ is not a square, then \eqref{eq:2} has at most three
coprime positive integer solutions.
\end{corollary}

\begin{remarks-nonum}
(1) For $m$ even, $a=1$, $p=2$ and $b=2^{m}$, with $X=2^{m/2+1}$ and $Y=2$, we have
$X^{2}- \left( a^{2}+b \right) Y^{2}=-4$. Also $\gcd \left( a^{2}, b \right)=1$,
$a^{2}+b$ is not a square and although $a^{2}+b \equiv 1 \pmod{8}$ if $m \geq 3$,
we have $b=2^{m}$ here. So the conditions in Corollary~\ref{cor:1.2-eqn1}
are satisfied in this case. Hence the results here both generalise and improve
the results in \cite{SWY}.

\vspace*{1.0mm}

(2) In Subsections~\ref{subsect:thm12-pm-egs}--\ref{subsect:thm12-16pm-egs},
we give examples showing that
Theorem~\ref{thm:1.2-seq} and Corollary~\ref{cor:1.2-eqn1} are best possible.

\vspace*{1.0mm}

(3) We have excluded $b=8p^{m}$ and $16p^{m}$ with $p$ odd from Corollary~\ref{cor:1.2-eqn1},
although they were included in Theorem~\ref{thm:1.2-seq}.
For $b=8p^{m}$ and $16p^{m}$ with $p$ odd and $a$ odd, we have $a^{2}+b \equiv 1 \pmod{8}$,
which is excluded in Lemma~\ref{lem:pell-rep}. While if $a$ is even, then
$4|\gcd \left ( a^{2}, b \right)$, so the condition that $\gcd \left( d_{2}, e \right)=1$
in Lemma~\ref{lem:pell-rep} is violated.

\vspace*{1.0mm}

(4) When $b=2p^{m}$ with $p$ odd, $X^{2}- \left( a^{2}+b \right) Y^{2}=-4$ is not
solvable when $a$ is odd, arguing modulo $8$. So $a$ must be even.
\end{remarks-nonum}

For small $m$, we have been able to remove the coprimality condition on the solutions
in Corollary~\ref{cor:1.2-eqn1}. We have proven the following.

\begin{corollary}
\label{cor:1.2-eqn2}
Let $a$, $m$ and $p$ be positive integers with $m=1,2,3,4$, and $p$ a prime. Put
$b=p^{m}$. Suppose that $\gcd \left( a^{2}, b \right)$ is squarefree,
$a^{2}+b$ is not a square and that $X^{2}-\left( a^{2}+b \right)Y^{2}=-4$
has a solution in positive integers.
Then \eqref{eq:2} has at most three positive integer solutions.
\end{corollary}

\begin{remarks-nonum}
(1) Corollary~\ref{cor:1.2-eqn2} is also best possible. See Subsection~\ref{subsect:cor12-egs}
for examples.

(2) The assumption that $a^{2}+b$ is not a square is not strictly necessary in either
Corollary~\ref{cor:1.2-eqn1} or \ref{cor:1.2-eqn2}. When $a^{2}+b$ is a positive
square, $X^{2}-\left( a^{2}+b \right)Y^{2}=\pm 4$ has no solution in positive
integers, as no two positive squares have difference $4$.
\end{remarks-nonum}

This paper uses the technique presented in \cite{V5}.
The proof there has its basis in the work of Siegel \cite{Siegel2} (also see
\cite{Evert1}). We shall refer to \cite{V5} and its lemmas throughout this
paper.

We structure this
paper as follows. In Section~\ref{sect:prelim}, we collect various facts that we
will require about squares in the sequence of $y_{k}$'s defined by \eqref{eq:yk-defn}.
In particular, we establish sharper lower bounds for $\left| y_{k} \right|$ for $k \neq 0$ and
an improved gap principle that hold under the conditions here (compare to
Lemmas~\lbLemma{} and \gapLemma{} in \cite{V5}).
In Section~\ref{sect:prop-11}, we bring together these results with the diophantine
approximation results in Section~\dioSect{} of \cite{V5}. As in Section~\genPropSection{}
of \cite{V5}, we state and prove Proposition~\ref{prop:4.1} here.
Our theorem and its corollaries follow from this proposition. Their
proofs are given in Section~\ref{sect:proofs}. Finally, in Section~\ref{sect:egs},
we provide examples showing the sharpness of our results.

Also, encouraged by the referee, PARI/GP code \cite{Pari} used for this work is
publicly available at \url{https://github.com/PV-314/hypgeom}. The author is
very happy to help interested readers who have any questions, problems or
suggestions for the use of this code.

\section{Lemmas about $\left( x_{k} \right)_{k=-\infty}^{\infty}$ and $\left( y_{k} \right)_{k=-\infty}^{\infty}$}
\label{sect:prelim}

For any non-zero integer, $n$, let $\core(n)$ be the unique squarefree integer
such that $n/\core(n)$ is a square. We will put $\core(1)=1$.

We start with a lower bound for $\left| y_{k} \right|$ when $k \neq 0$. For the
values of $b$ considered here, we get an improved dependence on $b$ in the
lower bounds compared to Lemma~\lbLemma{} of \cite{V5}: $O \left( b^{2} \right)$ here
versus $O(b)$ in Lemma~\lbLemma{} of \cite{V5}. These improved lower bounds are crucial
for enabling us to eliminate the possibility of squares among small values of
$y_{k} \neq 1$.

\begin{lemma}
\label{lem:Y-LB}
Suppose that $a$ and $b$ are positive integers with $a^{2}+b$ not a square
and let the $y_{k}$'s be defined by $\eqref{eq:yk-defn}$.

\noindent
{\rm (a)} Suppose that $b=p^{m}$ where $m$ is a non-negative integer, $p$ is a
prime and that $y_{k}$ is a square where $k \neq 0$. Then
\[
y_{k} \geq
\left\{
\begin{array}{ll}
b^{2}/16  & \text{if $p \neq 2$,}\\
b^{2}/2^{m+2}  & \text{if $b=2^{m} \leq 32$,}\\
b^{2}/256 & \text{if $b=2^{m} \geq 64$.}
\end{array}
\right.
\]

\noindent
{\rm (b)} Suppose that $b=2^{\ell}p^{m}$ with $\ell>0$, $m \geq 0$ and $p$ an
odd prime.
If $k \neq 0$ and $y_{k}$ is a square, then $y_{k}>b^{2}/64$ for $\ell=2$ and
$y_{k}>b^{2}/2^{2\ell}$ otherwise.
\end{lemma}

\begin{remark-nonum}
The results in each part of this lemma are best possible or else near to best
possible. See Subsection~\ref{subsect:YLB-egs} for examples demonstrating this.
\end{remark-nonum}

We start with a preliminary lemma that we will use to prove Lemma~\ref{lem:Y-LB}.

\begin{lemma}
\label{lem:Y-LB-1}
Let $a$, $b$, $t$, $u$ and $\varepsilon$ be as at the start of
Subsection~$\ref{subsect:results}$.
Suppose for some non-zero integer, $k$, $y_{k}$ is a perfect square. Then
\begin{equation}
\label{eq:g1-defn}
g_{1}= \gcd \left( \left( t_{k}+au_{k} \right)^{2}, bu_{k}^{2} \right)=1 \text{ or } 4,
\end{equation}
where $\varepsilon^{k}= \left( t_{k}+u_{k} \sqrt{a^{2}+b} \right)/2$ for integers
$t_{k}$ and $u_{k}$.
\end{lemma}

\begin{proof}
Putting $y_{k}=y^{2}$, and expanding the expression for $x_{k}+y_{k}\sqrt{a^{2}+b}$
in \eqref{eq:yk-defn}, we find that
\begin{equation}
\label{eq:15-a}
(2y)^{2} = \left( t_{k} + au_{k} \right)^{2} + bu_{k}^{2}.
\end{equation}

Since $\left( t_{k}+au_{k} \right)\left( t_{k}-au_{k} \right)=bu_{k}^{2} \pm 4$,
we know that
\[
\gcd \left( \left( t_{k}+au_{k} \right)\left( t_{k}-au_{k} \right), bu_{k}^{2} \right)|4.
\]

We now consider the various cases that can arise.

\vspace*{1.0mm}

(i) Suppose that $\gcd \left( \left( t_{k}+au_{k} \right)\left( t_{k}-au_{k} \right), bu_{k}^{2} \right)=1$.

In this case, one of the two quantities in the gcd must be odd. From
$\left( t_{k}+au_{k} \right)\left( t_{k}-au_{k} \right)=bu_{k}^{2} \pm 4$, they
must then both be odd. Hence every pair of $(2y)^{2}$, $\left( t_{k}+au_{k} \right)^{2}$
and $bu_{k}^{2}$ is relatively prime.

\vspace*{2.0mm}

(ii) Suppose that $\gcd \left( \left( t_{k}+au_{k} \right)\left( t_{k}-au_{k} \right), bu_{k}^{2} \right)=2$.

If $\gcd \left( \left( t_{k}+au_{k} \right)\left( t_{k}-au_{k} \right), bu_{k}^{2} \right)=2$,
then $t_{k}+au_{k}$ and $t_{k}-au_{k}$ are both even, so $2 \parallel \left( bu_{k}^{2} \right)$.
But this is impossible by $(2y)^{2} = \left( t_{k}+au_{k} \right)^{2} + bu_{k}^{2}$.

\vspace*{2.0mm}

(iii) Suppose that $\gcd \left( \left( t_{k}+au_{k} \right)\left( t_{k}-au_{k} \right), bu_{k}^{2} \right)=4$.

We consider separately the two cases of $au_{k}$ odd or even.

(iii-a) Suppose that $au_{k}$ is odd.

Then $t_{k}$ is odd and $4$ divides exactly one of $t_{k}+au_{k}$ and $t_{k}-au_{k}$.
This means that $4 \parallel \left( bu_{k}^{2} \right)$.

(iii-a-1) Suppose that $4 | \left( t_{k}+au_{k} \right)$.

Then $y$ is odd, as otherwise $16| bu_{k}^{2}$. Thus any pair of the three
quantities, $(2y)^{2}$, $\left( t_{k}+au_{k} \right)^{2}$ and $bu_{k}^{2}$
has gcd $4$.

(iii-a-2) Suppose that $2 \parallel \left( t_{k}+au_{k} \right)$.

Again, any
pair of the three quantities, $(2y)^{2}$, $\left( t_{k}+au_{k} \right)^{2}$ and
$bu_{k}^{2}$ has gcd $4$.

(iii-b) Suppose that $au_{k}$ is even.

Then $t_{k}+au_{k} \equiv t_{k}-au_{k} \pmod{4}$. So either $2 \parallel \left( t_{k}+au_{k} \right)$
and $2 \parallel \left( t_{k}-au_{k} \right)$ or else $16| \left( \left( t_{k}+au_{k} \right) \left( t_{k}-au_{k} \right) \right)$
and $4 \parallel \left( bu_{k}^{2} \right)$.

(iii-b-1) Suppose that $2 \parallel \left( t_{k}-au_{k} \right)$ and
$2 \parallel \left( t_{k}+au_{k} \right)$.

In this case, $4 \parallel \left( t_{k}+au_{k} \right)^{2}$, so any pair of the
three quantities, $(2y)^{2}$, $\left( t_{k}+au_{k} \right)^{2}$ and $bu_{k}^{2}$
has gcd $4$.

(iii-b-2) Suppose that $16| \left( \left( t_{k}+au_{k} \right) \left( t_{k}-au_{k} \right) \right)$
and $4 \parallel \left( bu_{k}^{2} \right)$.

In this case, since $4\parallel \left( bu_{k}^{2} \right)$, we saw above that any
pair of the three quantities, $(2y)^{2}$, $\left( t_{k}+au_{k} \right)^{2}$ and
$bu_{k}^{2}$ has gcd $4$.

\vspace*{1.0mm}

Hence from our consideration of cases~(i), (ii) and (iii) above, we have shown
that any pair of the three terms, $(2y)^{2}$, $\left( t_{k}+au_{k} \right)^{2}$
and $bu_{k}^{2}$ has gcd $1$ or $4$. Since $4|(2y)^{2}$, we need only be concerned
about $g_{1}$ and we have shown it can only take the values $1$ or $4$.
\end{proof}

\begin{proof}[Proof of Lemma~\ref{lem:Y-LB}]
We will use the same notation here as in Lemma~\ref{lem:Y-LB-1} and its proof.
From Lemma~\ref{lem:Y-LB-1}, we know that
\[
g_{1}= \gcd \left( \left( t_{k}+au_{k} \right)^{2}, bu_{k}^{2} \right)=1 \text{ or } 4.
\]

We consider both possibilities separately.

\vspace*{1.0mm}

(i) Suppose that $g_{1}=1$.

We can use this, along with the expression for $(2y)^{2}$ in \eqref{eq:15-a} to
see that there exist integers $b_{1}>0, b_{2}>0$,
$r$ and $s$ such that
\begin{equation}
\label{eq:y-g1}
2y - \left( t_{k}+au_{k} \right)=b_{1}r^{2}, \quad
2y + \left( t_{k}+au_{k} \right)=b_{2}s^{2},
\end{equation}
where $\gcd \left( b_{1}r^{2}, b_{2}s^{2} \right)=1$ and $b_{1}b_{2}=b$ is odd.

\vspace*{1.0mm}

(ii) Suppose that $g_{1}=4$.

We put
\[
g_{2}=\gcd \left( 2y-\left( t_{k}+au_{k} \right), 2y+\left( t_{k}+au_{k} \right) \right).
\]

From $g_{1}=4$ and the expression for $(2y)^{2}$ in \eqref{eq:15-a}, it follows that
$g_{2}=2$ or $4$ (if it were $8$, then $16$ would divide $g_{1}$).

There are three possibilities.

\vspace*{1.0mm}

(ii-1) Suppose that $g_{2}=4$.

In this case, there exist integers $b_{1}, b_{2}$, $r$ and $s$ such that
\begin{equation}
\label{eq:y-g4a}
2y - \left( t_{k}+au_{k} \right)=b_{1}r^{2}, \quad
2y + \left( t_{k}+au_{k} \right)=b_{2}s^{2},
\end{equation}
where $g_{2}=\gcd \left( b_{1}r^{2}, b_{2}s^{2} \right)=4$ and $b_{1}b_{2}=b$,
with either
$\gcd \left( r^{2}, s^{2} \right)=4$ and $\gcd \left( b_{1}, b_{2} \right)=1$,
or else
$\gcd \left( r^{2}, s^{2} \right)=1$ and $\gcd \left( b_{1}, b_{2} \right)=1,2,4$.

\vspace*{1.0mm}

(ii-2) Suppose that $g_{2}=2$ and $b$ is odd.

In this case, there exist integers $b_{1}, b_{2}$, $r$ and $s$ such that
\begin{equation}
\label{eq:y-g4b}
2y - \left( t_{k}+au_{k} \right)=2b_{1}r^{2}, \quad
2y + \left( t_{k}+au_{k} \right)=2b_{2}s^{2},
\end{equation}
where $\gcd \left( b_{1}r^{2}, b_{2}s^{2} \right)=1$ and $b_{1}b_{2}=b$.

\vspace*{1.0mm}

(ii-3) Suppose that $g_{2}=2$ and $b$ is even.

In this case, there exist integers $b_{1}, b_{2}$, $r$ and $s$ such that
\begin{equation}
\label{eq:y-g4c}
2y - \left( t_{k}+au_{k} \right)=b_{1}r^{2}, \quad
2y + \left( t_{k}+au_{k} \right)=b_{2}s^{2},
\end{equation}
where $\gcd \left( b_{1}r^{2}, b_{2}s^{2} \right)=2$ and $b_{1}b_{2}=b$.

Since $g_{2}=2$ and $2y$ is even, $t_{k}+au_{k}$ is also even. One of $2y$ or
$t_{k}+au_{k}$ must be $2 \bmod{4}$, while the other is $0 \bmod{4}$ (otherwise $g_{2} \neq 2$).
Hence, from the expression for $(2y)^{2}$ in \eqref{eq:15-a}, we see that $bu_{k}^{2} \equiv \pm 4 \pmod{16}$.
Since $b$ is even, it follows that $u_{k}$ is odd (otherwise $8| \left( bu_{k}^{2} \right)$).
So $2y \pm \left( t_{k}+au_{k} \right) \equiv 2 \pmod{4}$.
Since $u_{k}$ is odd, it follows that $\gcd \left( b_{1}, b_{2} \right)=2$.

\vspace*{1.0mm}

To prove both parts of this lemma, we consider each of the four possibilities
in equations~\eqref{eq:y-g1}--\eqref{eq:y-g4c}. Recall that we are only considering
$b=p^{m}$ for part~(a) and $b=2^{\ell}p^{m}$ for part~(b).

\vspace*{1.0mm}

(a-1) From \eqref{eq:y-g1} with $y=\sqrt{y_{k}}$, we have
$4\sqrt{y_{k}}=b_{1}r^{2}+b_{2}s^{2}$ with $r, s \geq 1$. Since $\gcd \left( b_{1}, b_{2} \right)=1$,
it follows that $b_{1}=b$ and $b_{2}=1$ (note that $b_{1}=1$ and $b_{2}=b$
is also possible, but without loss of generality, we need not consider this
possibility either here or in what follows). Thus $4\sqrt{y_{k}}>b$.

\vspace*{1.0mm}

(a-2) From \eqref{eq:y-g4a}, we again have $4\sqrt{y_{k}}=b_{1}r^{2}+b_{2}s^{2}$.
There are four subcases.

\vspace*{1.0mm}

(a-2-i) One possibility is that $\gcd(r,s)=2$ and
$\gcd \left( b_{1}, b_{2} \right)=1$, so $b_{1}=p^{m}$ and $b_{2}=1$,
$r, s \geq 2$, and hence, $4\sqrt{y_{k}}>4b$.

\vspace*{1.0mm}

(a-2-ii) If $\gcd(r,s)=1$ and $\gcd \left( b_{1}, b_{2} \right)=1$, then
$b_{1}=p^{m}$ and $b_{2}=1$, so $4\sqrt{y_{k}}>b$. This is consistent
with the inequalities in part~(a) except for $b=2$.

For $b=2$, we have $4\sqrt{y_{k}}=b_{1}r^{2}+b_{2}s^{2} \geq 2+1=3$.
So $y_{k} \geq 1=b^{2}/4$.

\vspace*{1.0mm}

(a-2-iii) If $\gcd(r,s)=1$ and $\gcd \left( b_{1}, b_{2} \right)=2$, then
$p=2$ and $m \geq 2$, so $b_{1}=2$ and $b_{2}=b/2$. Hence
$4\sqrt{y_{k}}>b/2$. This is consistent
with the inequalities in part~(a) except for $b=4$ or $b=8$.

For $b=4$, we have $b_{1}=b_{2}=2$, so $4\sqrt{y_{k}} \geq 2+2=b$.
Hence $y_{k} \geq b^{2}/16$.

For $b=8$, we have $b_{1}=4$ and $b_{2}=2$, so
$4\sqrt{y_{k}} \geq 2+4=3b/4$. Hence $y_{k} \geq 9b^{2}/256>b^{2}/32$.

\vspace*{1.0mm}

(a-2-iv) If $\gcd(r,s)=1$ and $\gcd \left( b_{1}, b_{2} \right)=4$, then
$p=2$ and $m \geq 4$, so $b_{1}=4$ and $b_{2}=b/4$. Hence
$4\sqrt{y_{k}}>b/4$. This is consistent with the
inequalities in part~(a) except for $b=16$ or $32$.

For $b=16$, we have $b_{1}=b_{2}=4$, so $4\sqrt{y_{k}} \geq 4+4=8=b/2$.
Hence $y_{k} \geq b^{2}/2^{m+2}$.

Similarly, for $b=32$, we have $b_{1}=4$ and $b_{2}=8$ (or vice versa), so
$4\sqrt{y_{k}} \geq 4+8 =12=12b/32$. Hence $y_{k} \geq 9b^{2}/1024>b^{2}/2^{m+2}$.

\vspace*{1.0mm}

(a-3) From \eqref{eq:y-g4b}, we have $4\sqrt{y_{k}}=2b_{1}r^{2}+2b_{2}s^{2}$
with $r, s \geq 1$, $b_{1}=p^{m}$ and $b_{2}=1$. Thus $4\sqrt{y_{k}}>2b$.

\vspace*{1.0mm}

(a-4) From \eqref{eq:y-g4c}, we have $4\sqrt{y_{k}}=b_{1}r^{2}+b_{2}s^{2}$ with
$r, s \geq 1$, $b_{1}=2^{m-1}$ and $b_{2}=2$. Thus $4\sqrt{y_{k}}>b/2$.
This is consistent with the inequalities in part~(a) except for $b=4,8$.

The treatment of $b=4$ and $8$ is as in case~(a-2-iii) above.

Thus the inequality in part~(a) always holds.

\vspace*{3.0mm}

(b) The proof is very similar to the proof of part~(a). Since $b$ is even, we
saw above that we need only consider \eqref{eq:y-g4a}
and \eqref{eq:y-g4c}.

\vspace*{1.0mm}

(b-1) From \eqref{eq:y-g4a}, we have $4\sqrt{y_{k}}=b_{1}r^{2}+b_{2}s^{2}$ with either
$r, s \geq 2$ (when $\gcd(r,s)=2$ and $\gcd \left( b_{1}, b_{2} \right)=1$) or
else $r,s \geq 1$ (when $\gcd(r,s)=1$ and $\gcd \left( b_{1}, b_{2} \right)=1,2,4$).

\vspace*{1.0mm}

(b-1-i) In the case that $r, s \geq 2$ (when $\gcd(r,s)=2$ and $\gcd \left( b_{1}, b_{2} \right)=1$),
since $b=2^{\ell}p^{m}$, it follows that $b_{1}=p^{m}$ and $b_{2}=2^{\ell}$ or
$b_{1}=1$ and $b_{2}=b$.
If $b_{1}=p^{m}$ and $b_{2}=2^{\ell}$, then
$4\sqrt{y_{k}}>4\max \left( 2^{\ell}, p^{m} \right)$
and so $y_{k}>b^{2}/2^{2\ell}$.
If $b_{1}=1$ and $b_{2}=b$, then
$4\sqrt{y_{k}}>4b$ and so $y_{k}>b^{2}$.

\vspace*{1.0mm}

(b-1-ii) If $r,s \geq 1$ and $\gcd \left( b_{1}, b_{2} \right)=1$, then there are
two possibilities.

The first possibility is that $b_{1}=b$ and $b_{2}=1$, so $4\sqrt{y_{k}}>b$. Thus
$y_{k}>b^{2}/16$.
This is consistent with the inequalities in part~(b) except for $\ell=1$.

For $\ell=1$, from $4\sqrt{y_{k}}=br^{2}+s^{2}$, we observe that $s$ must be
even. Therefore $4| \left( br^{2} \right)$, which implies that $r$ is even.
Hence $\sqrt{y_{k}}>b$, which is consistent with the inequalities in part~(b).

The second possibility is that $b_{1}=p^{m}$ and $b_{2}=2^{\ell}$. Since $g_{1}=4$,
both $2y \pm \left( t_{k}+au_{k} \right)$ are even. Combining this with the
expression for $b_{1}$, we see that $r$ must be even.
So $4\sqrt{y_{k}} \geq 4b_{1}+b_{2}>4p^{m}$.
and $y_{k}>b^{2}/2^{2\ell}$.

\vspace*{1.0mm}

(b-1-iii) If $r,s \geq 1$ and $\gcd \left( b_{1}, b_{2} \right)=2$, then there are
two possibilities.

First, we could have
$b_{1}=2^{\ell-1}p^{m}$ and $b_{2}=2$ with $\ell \geq 2$,
so $4\sqrt{y_{k}}>b/2$ (i.e., $y_{k}>b^{2}/64$, which is consistent with the
inequalities in part~(b)).

Second, we could have $b_{1}=2p^{m}$ and
$b_{2}=2^{\ell-1}$ with $\ell \geq 2$,
so $4\sqrt{y_{k}}>\max \left( 2^{\ell-1}, 2p^{m} \right)$
and $y_{k}>b^{2}/2^{2\ell+2}$.
For $\ell=2$, this is consistent with the inequalities in part~(b).
If $\ell \geq 3$, then from $4\sqrt{y_{k}}=b_{1}r^{2}+b_{2}s^{2}$, we see that
$4 | \left( b_{1}r^{2} \right)$. So $r$ must be even. Therefore,
$y_{k}>4 b^{2}/2^{2\ell}$.

\vspace*{1.0mm}

(b-1-iv) If $r,s \geq 1$ and $\gcd \left( b_{1}, b_{2} \right)=4$, then there are
two possibilities.

We could have
$b_{1}=2^{\ell-2}p^{m}$ and $b_{2}=4$ with $\ell \geq 4$,
so $4\sqrt{y_{k}}>2^{\ell-2}p^{m}$ (i.e., $y_{k}>b^{2}/256 \geq b^{2}/2^{2\ell}$),

Or we could have $b_{1}=4p^{m}$ and
$b_{2}=2^{\ell-2}$,
so $4\sqrt{y_{k}}>\max \left( 2^{\ell-2}, 4p^{m} \right)$
and $y_{k}>b^{2}/2^{2\ell}$.

\vspace*{1.0mm}

(b-2) From \eqref{eq:y-g4c}, we have $4\sqrt{y_{k}}=b_{1}r^{2}+b_{2}s^{2}$ with
$r, s \geq 1$ and $\gcd \left( b_{1}, b_{2} \right)=2$. We must have
$b=2^{\ell}p^{m}$ with $\ell \geq 2$, $b_{1}=2^{\ell-1}$ and $b_{2}=2p^{m}$
or else $b_{1}=2$ and $b_{2}=2^{\ell-1}p^{m}$.

\vspace*{1.0mm}

(b-2-i) In the first case, if $\ell \geq 3$,
then $2y - \left( t_{k}+au_{k} \right)=b_{1}r^{2} \equiv 0 \pmod{4}$.
Since $g_{2}= \gcd \left( 2y - \left( t_{k}+au_{k} \right), 2y + \left( t_{k}+au_{k} \right) \right)=2$,
it follows that $2y + \left( t_{k}+au_{k} \right)=b_{2}s^{2} \equiv 2 \pmod{4}$ -- but $4|b_{2}$,
so this is not possible.
So this case can only occur for $\ell=2$.
Hence $4\sqrt{y_{k}}>\max \left( 2, 2p^{m} \right)$, and
$y_{k}>b^{2}/64$.

\vspace*{1.0mm}

(b-2-ii) In the second case,
$4\sqrt{y_{k}}>2^{\ell-1}p^{m}$, so $y_{k}>b^{2}/64$. Notice that for $\ell \geq 3$,
we have $b^{2}/64 \geq b^{2}/2^{2\ell}$.

Thus the inequality in part~(b) always holds.
\end{proof}

\begin{lemma}
\label{lem:gap}
Suppose that $a$ and $b$ are positive integers and $a^{2}+b$ is not a
square. Let the $y_{k}$'s be defined as in $\eqref{eq:yk-defn}$.

{\rm (a)} Suppose that $b$ is a square. If $y_{i}$ and $y_{j}$
are both squares with $y_{j}>y_{i}>1$, then
\[
y_{j} > 57.32 \left( \frac{a^{2}+b}{b} \right)^{2} y_{i}^{3}.
\]

{\rm (b)} Suppose $b=2^{\ell}p^{m}$ is not a square. If
$y_{k_{1}}$, $y_{k_{2}}$ and $y_{k_{3}}$ are three distinct squares with
$y_{k_{3}}>y_{k_{2}}>y_{k_{1}}>1$, then there exist distinct
$i,j \in \{ k_{1},k_{2}, k_{3} \}$ such that
\[
y_{j} > 15 \left( \frac{a^{2}+b}{b} \right)^{2} y_{i}^{3}.
\]
\end{lemma}

\begin{proof}
(a) This is Lemma~\gapLemma(a) of \cite{V5} upon recalling that $N_{\alpha}$
there is $-b$ here, $d$ there is $a^{2}+b$ here and that $b$ there is $1$ here.
The conditions on $y_{j}$ and $y_{i}$ Lemma~\gapLemma(a) of \cite{V5} hold
since here we must have $y_{j}>y_{i} \geq 4$.

(b) Our proof here will follow closely the proof of Lemma~\gapLemma(b) of \cite{V5}.
As in the proof of Lemma~\gapLemma{} of \cite[just above equation~(3.19) on page~314]{V5},
we put $\omega= \left( a-\sqrt{-b} \right) / \left( a+\sqrt{-b} \right)$.
We can write $\omega=e^{i \varphi}$, where $-\pi<\varphi \leq \pi$.
For any real number $\nu$, we shall put $\omega^{\nu}=e^{i\nu \varphi}$.

For any $j$ among $k_{1}$, $k_{2}$ and $k_{3}$, let $\zeta_{4}^{(j)}$ be the
$4$-th root of unity such that
\[
\left| \omega^{1/4} - \zeta_{4}^{(j)} \frac{r_{j}+s_{j}\sqrt{-\core(b)}}{r_{j}-s_{j}\sqrt{-\core(b)}} \right|
= \min_{0 \leq k \leq 3} \left| \omega^{1/4} - e^{2k \pi i/4} \frac{r_{j}+s_{j}\sqrt{-\core(b)}}{r_{j}-s_{j}\sqrt{-\core(b)}} \right|.
\]
where $r_{j}$ and $s_{j}$ are the values of $r$ and $s$ in Proposition~\repProp{} of \cite{V5}
associated with $x=x_{j}$ and $y=\sqrt{y_{j}}$.

Also for any such $j$, let $f_{j}$ be the value of $f$ in Proposition~\repProp{}
of \cite{V5} associated with $x=x_{j}$ and $y=\sqrt{y_{j}}$, where we use
$\alpha=a+\sqrt{a^{2}+b}$, so $N_{\alpha}=-b$. From Proposition~\repProp{}(c) of
\cite{V5}, we have $f_{j}|4$.

Choose $i$ and $j$ from among $k_{1}$, $k_{2}$ and $k_{3}$ so that $\zeta_{4}^{(i)}= \pm \zeta_{4}^{(j)}$.

Suppose $f_{i}f_{j} \leq 4$. From Lemma~\gapLemma(b) of \cite{V5}, we have
\[
y_{j}>15.36 \left( \frac{a^{2}+b}{4b} \right)^{2} y_{i}^{3}.
\]

Since $y_{i} \geq 4$, we obtain $y_{j}>15.36y_{i}$. Equation~(\gpEqnb) of \cite{V5}
applied with our notation and assumptions here states that
\begin{equation}
\label{eq:gpEqnb}
\frac{2}{\sqrt{f_{i}f_{j}} \left( y_{i}y_{j} \right)^{1/4}}
< 0.5051 \frac{\sqrt{b}}{\sqrt{a^{2}+b}}
\left( \frac{1}{y_{i}} + \frac{1}{y_{j}} \right).
\end{equation}

Applying $f_{i}f_{j} \leq 4$ to this inequality and $y_{j}>15.36y_{i}$ to
the right-hand side of it, it follows that
\begin{equation}
\label{eq:yj-LB2}
y_{j}>11.93 \left( \frac{a^{2}+b}{b} \right)^{2} y_{i}^{3}.
\end{equation}

From \eqref{eq:yj-LB2} and $y_{i} \geq 4$, we have $y_{j}>190.88y_{i}$. Applying
this to \eqref{eq:gpEqnb}, as we did above with $y_{j}>15.36y_{i}$, we obtain
\[
y_{j}>15.04 \left( \frac{a^{2}+b}{b} \right)^{2} y_{i}^{3}.
\]

Now suppose that $f_{i}f_{j} \geq 8$. We show that we are able to improve
\eqref{eq:gpEqnb} here. In this case, both $f_{i}$ and $f_{j}$
must be even and either $f_{i}=4$ or $f_{j}=4$. From Proposition~\repProp(c)
of \cite{V5}, this can only occur when $b \equiv 3 \pmod{4}$ and $4| \left( a^{2}+b \right)$
(i.e., $a$ is odd).

From the expression for $fy$ in equation~(3.2) in Proposition~\repProp{} of \cite{V5},
we see that here $r_{i}$ and $s_{i}$, as well as $r_{j}$ and $s_{j}$, must have
the same parity (since $f_{i}$ and $f_{j}$ are both even, $b$ in the notation of
\cite{V5} is $1$ and $b$ in the present paper satisfies $b \equiv 3 \pmod{4}$, so
$\core(b)$ odd). In the notation in the proof of Lemma~\gapLemma{} of \cite{V5}
(see immediately above equation~(\eqZetaIJRel{}) there), we write
\begin{align*}
x+y\sqrt{\core(-b)}
&= \left( r_{i}+s_{i}\sqrt{\core(-b)} \right) \left( r_{j}+s_{j}\sqrt{\core(-b)} \right) \\
&=r_{i}r_{j}+s_{i}s_{j}\core(-b) + \left( r_{i}s_{j}+r_{j}s_{i} \right) \sqrt{\core(-b)}.
\end{align*}

Since $r_{i}$ and $s_{i}$ have the same parity, as do $r_{j}$ and $s_{j}$, and
since $N_{\alpha}$ is odd, we see that $x$ and $y$ are both even here.
Therefore, the expressions for the numerators in equation~(\gpNumEqn) in \cite{V5}
have absolute values at least $4$ when $\zeta_{4}^{(j)}=\pm \zeta_{4}^{(i)}$ and
so
\[
\frac{4}{\sqrt{f_{i}f_{j}} \left( y_{i}y_{j} \right)^{1/4}}
< 0.5051 \frac{\sqrt{b}}{\sqrt{a^{2}+b}}
\left( \frac{1}{y_{i}} + \frac{1}{y_{j}} \right),
\]

The worst case is $f_{i}=f_{j}=4$. Here, since $y_{j}>y_{i}$, the
right-hand side is less than
$0.5051 \sqrt{b/ \left( a^{2}+b \right)} \left( 2/y_{i} \right)$,
and we obtain $y_{j}>0.96 \left( \left( a^{2}+b \right)/b \right)^{2}y_{i}^{3}$.
Since $y_{i} \geq 4$, this yields $y_{j}>15.36y_{i}$.

Using the same process as for $f_{i}f_{j} \leq 4$
(except that we need three passes through that process, rather than two passes
as we used for $f_{i}f_{j} \leq 4$), we obtain our result.
\end{proof}

It is the following lemma that allows us to associate solutions of the diophantine
equations \eqref{eq:2} (specifically under the conditions in Corollary~\ref{cor:1.2-eqn1})
with squares in these sequences.

\begin{lemma}
\label{lem:pell-rep}
Let $d \geq 2$, $m_{1}, m_{2} \geq 0$ and $p \geq 3$ be integers with $p$ a prime,
$d=d_{1}d_{2}^{2}$, where $d_{1}$ is squarefree, $d$ is not a square and
$d_{1} \not\equiv 1 \pmod{8}$ if both $m_{1} \geq 1$ and $m_{2} \geq 1$. Put
$e= \pm 2^{m_{1}}p^{m_{2}}$ and assume that $\gcd \left( d_{2}, e \right)=1$.

If $X^{2}-dY^{2}=e$ has a coprime positive integer solution $\left( x_{0}, y_{0} \right)$,
then all coprime positive integer solutions $(x,y)$ of $X^{2}-dY^{2}=e$ satisfy
\begin{equation}
\label{eq:14d}
\pm x + y \sqrt{d} = \left( x_{0}+y_{0}\sqrt{d} \right) \eta,
\end{equation}
for some choice of sign, where $\eta$ is a unit of norm $1$ in $\cO_{\bbQ \left( \sqrt{d} \right)}$.

If $e<0$, then we can choose $\eta>0$ and hence $\eta$ is a power
of the fundamental unit of $\cO_{\bbQ(\sqrt{d})}$.

Furthermore, we can write $\eta=\left( t+u\sqrt{d} \right)/2$ with $t,u \in \bbZ$.
\end{lemma}

\begin{remarks-nonum}
(1) The form of \eqref{eq:14d}, in particular the term $+y\sqrt{d}$ on the
left-hand side, is important to us for ensuring that solutions of the Diophantine
equations we consider in Corollary~\ref{cor:1.2-eqn1} come from squares in one
single sequence of $y_{k}$'s. More commonly, one would consider $x \pm y \sqrt{d}$
on the left-hand side of \eqref{eq:14d}.

\vspace*{2.0mm}

(2) The condition $d_{1} \not\equiv 1 \pmod{8}$ is required here, as otherwise the
ideal $(2)$ splits completely in $\bbQ \left( \sqrt{d} \right)$.
This can give rise to other families of solutions of $X^{2}-dY^{2}=e$.
For example, with $e=104$ and $d=d_{1}=17$, the solutions $(x,y)=(23,5)$ and
$\left( x_{0},y_{0} \right)=(11,1)$ do not satisfy the relationship in the lemma.

\vspace*{2.0mm}

(3) Similarly, the condition that $\gcd \left( d_{2},e \right)=1$
is required here too. For example, with $e=-185488$ and $d=867764$, where $d_{2}=2$,
the solutions $(x,y)=(32218340966, 34586161)$ and $\left( x_{0}, y_{0} \right)=(826,1)$
are related via
\[
x + y \sqrt{d} = \left( x_{0}+y_{0}\sqrt{d} \right) \left( \frac{36663027+78715\sqrt{d_{1}}}{2} \right),
\]
rather than via a unit of the form $\left( t+u\sqrt{d} \right)/2$.
\end{remarks-nonum}

\begin{proof}
This result is part of the classical literature for $e=\pm 1$.
So we will assume that at least one of $m_{1}$ or $m_{2}$ is positive in what
follows.

The proof uses the fact that for $\beta, \gamma \in \cO_{\bbK}$ for any number
field, $\bbK$, we have $(\beta)=(\gamma)$ if and only if $\beta=\gamma\varepsilon$
where $\varepsilon$ is a unit in $\cO_{\bbK}$.

In what follows, we let
$\bbK = \bbQ \left( \sqrt{d} \right)$ and, for any ideal, $\fa$, in $\cO_{\bbK}$,
$\overline{\fa}=\left\{ \sigma(\gamma): \gamma \in \fa \right\}$, where $\sigma$
is the non-identity element of the Galois group of $\bbK$.

Let $(x,y)$ be any relatively prime solution of $X^{2}-dY^{2}=e$. We will consider
$x+y \sqrt{d}$ and the principal ideal it generates, $\left( x+y \sqrt{d} \right)$.

We break our proof of \eqref{eq:14d} and the second statement of the lemma into
three steps, depending on the factorisation of $e$. At the end of the proof, we
establish the statement regarding the form of $\eta$.

\vspace*{1.0mm}

(i) Suppose that $e=\pm p^{m_{2}}$ with $m_{2} \geq 1$.

There are three possibilities for how the ideal $(p)$ factors into prime ideals
in $\cO_{\bbK}$. Either $(p)$ remains prime in $\cO_{\bbK}$, or
$(p)=\fp^{2}$ for a prime ideal $\fp$ in $\cO_{\bbK}$,
or else $(p)=\fp\overline{\fp}$, where $\fp$ is a prime ideal in $\cO_{\bbK}$ and
$\overline{\fp}$ is as defined above.

(i-1) If $(p)$ is a prime ideal in $\cO_{\bbK}$, from the ideal equation
$\left( x+y\sqrt{d} \right) \left( x-y\sqrt{d} \right) =(p)^{m_{2}}$,
$(p)=\overline{(p)}$ and the unique factorisation of ideals, it follows that
$\left( x+y\sqrt{d} \right)=(p)^{m_{2}/2}$. In the same way, we have
$\left( x_{0}+y_{0}\sqrt{d} \right)=(p)^{m_{2}/2}$. Hence
$\left( x+y\sqrt{d} \right)=\left( x_{0}+y_{0}\sqrt{d} \right)=(p)^{m_{2}/2}$.
So $m_{2}$ must be a positive even integer and $p$ must also divide the algebraic
integer, $x+y\sqrt{d}$.
Since $p \neq 2$ and $\gcd \left( p, d_{2} \right)=\gcd \left( e, d_{2} \right)=1$,
this implies that $p|x$ and $p|y$. But this is not possible, since $\gcd(x,y)=1$.

\vspace*{1.0mm}

(i-2) The second possibility is that $(p)=\fp^{2}$ (i.e., $p$ ramifies in $\cO_{\bbK}$).
Here $\fp=\overline{\fp}$, since $(p)=N(\fp)=\fp\overline{\fp}$ by the definition
of the norm of an ideal, and
then applying the unique factorisation of ideals. We use the ideal equation
$\left( x+y\sqrt{d} \right) \left( x-y\sqrt{d} \right) =(p)^{m_{2}}=\fp^{2m_{2}}$.
So $\left( x+y\sqrt{d} \right)=\fp^{m_{2}}$.
Similarly, $\left( x_{0}+y_{0}\sqrt{d} \right)=\fp^{m_{2}}$. Thus
$\left( x_{0}+y_{0}\sqrt{d} \right)=\left( x+y\sqrt{d} \right)$.
Therefore, by our observation at the start of the
proof, $x+y \sqrt{d}$ must be a unit times $x_{0}+y_{0}\sqrt{d}$. Since the
algebraic numbers $x+y \sqrt{d}$ and $x_{0}+y_{0}\sqrt{d}$ have the same norm
($\pm p^{m_{2}}$), this unit must be of norm $1$ and the relationship in
\eqref{eq:14d} follows.

Since
$x,y,x_{0}$ and $y_{0}$ are all positive integers, it follows from \eqref{eq:14d}
that $\eta$ must also be positive, so the second statement in the lemma also holds.

\vspace*{1.0mm}

(i-3) Lastly, we consider the case when $(p)$ splits in $\bbK$. We use the ideal equation
$\left( x+y\sqrt{d} \right) \left( x-y\sqrt{d} \right) =(p)^{m_{2}}=\fp^{m_{2}}
\overline{\fp}^{m_{2}}$.

Note that $x+y \sqrt{d}$ cannot be a
member of $\fp^{m} \overline{\fp}^{m_{2}-m}$ for $1 \leq m<m_{2}$, as such an ideal
would have a power of $(p)$ as a factor. As in case~(i-1), since $p \neq 2$,
$\gcd \left( p, d_{2} \right)=\gcd \left( e, d_{2} \right)=1$ and $\gcd(x,y)=1$,
this is not possible.

Hence we must have either $\left( x+y \sqrt{d} \right) = \fp^{m_{2}}$
or $\left( x+y \sqrt{d} \right) = \overline{\fp}^{m_{2}}$.
Similarly, $\left( x_{0}+y_{0}\sqrt{d} \right)=\fp^{m_{2}}$
or $\left( x_{0}+y_{0} \sqrt{d} \right) = \overline{\fp}^{m_{2}}$.
This implies that either
$\left( x+y\sqrt{d} \right)= \left( x_{0}+y_{0}\sqrt{d} \right)$ or
$\left( x+y\sqrt{d} \right)= \overline{\left( x_{0}+y_{0}\sqrt{d} \right)}
= \left( x_{0}-y_{0}\sqrt{d} \right)$.

We have already dealt with the first possibility when considering case~(i-2) above.
For the second possibility, we can write
$x+y\sqrt{d} = \left( x_{0}-y_{0}\sqrt{d} \right) \nu$, for some unit $\nu \in \cO_{\bbK}$.
Furthermore, as in case~(i-2), $\nu$ must be of norm $1$.
Taking conjugates and
multiplying both sides by $-1$, we have
$-x+y\sqrt{d} = \left( x_{0}+y_{0}\sqrt{d} \right) \left( -\overline{\nu} \right)
=\left( x_{0}+y_{0}\sqrt{d} \right) \left( -\nu^{-1} \right)$, since the norm of
$\nu$ is $1$. This establishes the relationship in \eqref{eq:14d} in this case.

Since $x,y,x_{0}$ and $y_{0}$ are all positive integers, if $e<0$,
then it follows that $x-y\sqrt{d}<0$. Hence $-x+y\sqrt{d}$ is positive, so, by
\eqref{eq:14d}, $\eta=-\nu^{-1}$ must be positive too. Hence the second statement
in the lemma also holds in this case.

\vspace*{1.0mm}

(ii) Suppose that $e=\pm 2^{m_{1}}$ with $m_{1} \geq 1$.

If $(2)$ is a prime ideal in $\cO_{\bbK}$ or if it ramifies in $\cO_{\bbK}$, then
we can argue as in cases~(i-1) and (i-2) above that
$\left( x+y \sqrt{d} \right) = \left( x_{0}+y_{0} \sqrt{d} \right)$.
The proof of the relationship in \eqref{eq:14d} and the second statement in the
lemma now follows as in case~(i-2). Note that, unlike in case~(i-1), the case
when $(2)$ is a prime ideal in $\cO_{\bbK}$ can occur. However, the treatment in
case~(i-2) still works in this case.

This leaves us with $(2)=\fq \overline{\fq}$.
This can only occur if $d_{1} \equiv 1 \pmod{8}$.
Since $\gcd \left( d_{2}, e \right)=1$ and $\gcd(x,y)=1$, it follows that
$x$ and $y$ must both be odd and that $\left( x+y\sqrt{d} \right)/2 \in \cO_{\bbK}$.
Hence $\left( x+y \sqrt{d} \right)$ must be $\fq \overline{\fq}$ times an ideal
in $\cO_{\bbK}$. That is,
$\left( x+y \sqrt{d} \right)$ must be of the form $\fq^{m} \overline{\fq}^{m_{1}-m}$
for some $1 \leq m < m_{1}$. If $m, m_{1}-m \geq 2$, then $(2)^{2}$ is a factor of
the ideal $\left( x+y \sqrt{d} \right)$, this is not possible since
$\gcd \left( 2, d_{2} \right)=\gcd(x,y)=1$. Hence
$\left( x+y \sqrt{d} \right) = \fq \overline{\fq}^{m_{1}-1}
=(2) \overline{\fq}^{m_{1}-2}$ or $\overline{\fq} \fq^{m_{1}-1}=(2) \fq^{m_{1}-2}$.
Since $x$ and $y$ are both odd, we have
$\left( \left( x+y \sqrt{d} \right)/2 \right) = \overline{\fq}^{m_{1}-2}$
or $\fq^{m_{1}-2}$. Similarly,
$\left( \left( x_{0}+y_{0} \sqrt{d} \right)/2 \right) = \overline{\fq}^{m_{1}-2}$
or $\fq^{m_{1}-2}$. Hence
$\left( \left( x+y \sqrt{d} \right)/2 \right)
=\left( \left( x_{0}+y_{0} \sqrt{d} \right)/2 \right)$ or
$\overline{\left( \left( x_{0}+y_{0} \sqrt{d} \right)/2 \right)}$.
Proceeding as in the proofs of cases~(i-2) and (i-3) completes the proof of
\eqref{eq:14d} and the second statement in the lemma in this case.

\vspace*{1.0mm}

(iii) Lastly, we suppose that $e=\pm 2^{m_{1}}p^{m_{2}}$ with $m_{1}, m_{2} \geq 1$.

Since $d_{1} \not\equiv 1 \pmod{8}$, $(2)$ is either a prime ideal in $\cO_{\bbK}$
or else the square of a prime ideal, $\fq$. Also, we saw in case~(i-1) that we
need not consider the case when $(p)$ is a prime ideal, so $(p)$
either ramifies in $\cO_{\bbK}$ or splits into the product of two prime ideals,
$\fp$ and $\overline{\fp}$. Since $\fp=\overline{\fp}$ when $(p)$ ramifies in
$\cO_{\bbK}$ (see case~(i-2) above), we can write $(p)=\fp \overline{\fp}$ in
both cases.

(iii-1) We first consider the case when $(2)$ is a prime ideal.

We start by showing that $N_{\bbK/\bbQ}\left( x+y \sqrt{d} \right)$ cannot be
$\pm 2p^{m_{2}}$. Since $(2)$ is a
prime ideal, it follows that $d_{1} \equiv 5 \pmod{8}$. So if $d$ is odd, then
$d \equiv 5 \pmod{8}$.
In this case, $x^{2}-dy^{2} \equiv 2 \pmod{4}$ is only possible if both $x$ and
$y$ are even, which we do not allow here. If $d$ is even, then $d_{2}$ must be
even, so $4|d$. Also $x$ must be
even, so $X^{2}-dY^{2}$ is divisible by $4$, eliminating this case.

Hence $N_{\bbK/\bbQ}\left( x+y \sqrt{d} \right)$ must be $\pm 2^{m_{1}}p^{m_{2}}$
with $m_{1} \geq 2$.
As in case~(i-3), $\left( x+y \sqrt{d} \right)$ cannot equal $(2)^{m_{1}/2}\fp^{m} \overline{\fp}^{m_{2}-m}$
for $1 \leq m<m_{2}$, unless $\fp=\overline{\fp}$. So it must be the case
that $\left( x+y \sqrt{d} \right) = (2)^{m_{1}/2}\fp^{m_{2}}$ or
$(2)^{m_{1}/2} \overline{\fp}^{m_{2}}$,
depending on how $(p)$ factors in $\cO_{\bbK}$.
Similarly,
$\left( x_{0}+y_{0} \sqrt{d} \right) = (2)^{m_{1}/2}\fp^{m_{2}}$ or
$(2)^{m_{1}/2} \overline{\fp}^{m_{2}}$. Hence either
$\left( x+y \sqrt{d} \right) = \left( x_{0}+y_{0} \sqrt{d} \right)$ or
$\overline{\left( x_{0}+y_{0} \sqrt{d} \right)}$.
Arguing as above in cases~(i-2) and (i-3), the relationship in \eqref{eq:14d}
and the second statement in the lemma hold in this case.

(iii-2) Now suppose that $(2)=\fq^{2}$, where $\fq$ is a prime ideal.

As in case~(i-2), we have $\fq=\overline{\fq}$ here. Since $(2)$ ramifies here
and $\gcd \left( d_{2}, e \right)=\gcd(x,y)=1$, it follows that
$N_{\bbK/\bbQ}\left( x+y \sqrt{d} \right)=e$
requires $e=\pm 2p^{m_{2}}$, rather than $e=\pm 2^{m_{1}}p^{m_{2}}$ with
$m_{1} \geq 2$. This holds as follows.

The ideal $(2)$ can only ramify in $\cO_{\bbK}$ if
$d_{1} \equiv 3 \pmod{4}$ or if $d_{1}$ is even. In the former case and if $d$ is
odd, then $x^{2}-dy^{2}
\equiv x^{2}+y^{2} \equiv 0 \pmod{4}$ can only happen if both $x$ and $y$ are
even, which we do not allow here.
We do not allow $d_{1} \equiv 3 \pmod{4}$ and $d$ is even, since $d_{2}$ is even
in this case and hence $\gcd \left( d_{2}, e \right) \neq 1$.

If $d_{1}$ is even, then $x^{2}-dy^{2} \equiv 0 \pmod{4}$, implies that $x$ is
also even. But if $d \equiv 2 \pmod{4}$ (i.e., $d_{2}$ is odd), then $y$ must
be even too, which we do not allow. Once again, we do not permit $d_{2}$ to be
even here.

As in case~(i-3), $\left( x+y \sqrt{d} \right)$ cannot equal $\fq\fp^{m} \overline{\fp}^{m_{2}-m}$
for $1 \leq m<m_{2}$, unless $\fp=\overline{\fp}$. So we must have either
$\left( x+y \sqrt{d} \right) = \fq\fp^{m_{2}}$ or $\left( x+y \sqrt{d} \right)
= \fq\overline{\fp}^{m_{2}} = \overline{\fq \fp^{m_{2}}}$, recalling that
$\fq=\overline{\fq}$ since $(2)$ ramifies. Similarly,
$\left( x_{0}+y_{0} \sqrt{d} \right) = \fq\fp^{m_{2}}$ or
$\overline{\fq \fp^{m_{2}}}$. Hence either
$\left( x+y \sqrt{d} \right) = \left( x_{0}+y_{0} \sqrt{d} \right)$ or
$\overline{\left( x_{0}+y_{0} \sqrt{d} \right)}$.
Arguing as above in cases~(i-2) and (i-3), the relationship in \eqref{eq:14d}
and the second statement in the lemma hold in this case.

\vspace*{1.0mm}

We now prove the last statement in the lemma.

We have shown that we can write
\[
x+y\sqrt{d} = \left( x_{0}+y_{0}\sqrt{d} \right) \left( \frac{t+u\sqrt{d_{1}}}{2} \right)
= \frac{x_{0}t+y_{0}ud_{1}d_{2}}{2} + \left( \frac{x_{0}u}{2d_{2}} + \frac{ty_{0}}{2} \right) \sqrt{d}.
\]

In order for $y$ to be a rational integer, we must have $d_{2} | \left( x_{0}u \right)$.
Since $\gcd \left( d_{2}, e \right)=1$ and $x_{0}^{2}-dy_{0}^{2}=e$, we must have
$\gcd \left( d_{2}, x_{0} \right)=1$ and so $d_{2}|u$. That is $u=u_{1}d_{2}$.
In this case, we can write $\left( t+u\sqrt{d_{1}} \right)/2=\left( t+u_{1}\sqrt{d} \right)/2$,
as required.
\end{proof}

\section{Proposition~\ref{prop:4.1} and its proof}
\label{sect:prop-11}

All of our results in Subsection~\ref{subsect:results} follow from
the next proposition.

\begin{proposition}
\label{prop:4.1}
Let $a$ and $b$ be positive integers such that $a^{2}+b$ is not
a square and let $\left( y_{k} \right)_{k=-\infty}^{\infty}$
be defined by \eqref{eq:yk-defn}.

\noindent
{\rm (a)} If $b$ is a square, then there is at most one integer
square among all distinct elements of $\left( y_{k} \right)_{k=-\infty}^{\infty}$
which satisfy $y_{k}>\max \left( 1,4.75b^{3/2}/\sqrt{a^{2}+b} \right)$ for $b$
even and $y_{k}>\max \left( 1,19b^{3/2}/\sqrt{a^{2}+b} \right)$ for $b$ odd.

\noindent
{\rm (b)} If $b=p^{m}$, where $p$ is an odd prime and $m$ is a positive
odd integer, then there are at most two distinct integer squares among all 
$y_{k}>\max \left( 1,b^{2}/16 \right)$.

\noindent
{\rm (c)} If $b$ is an even non-square of the form $2^{m}$, $2p^{m}$, $4p^{m}$,
$8p^{m}$ or $16p^{m}$, where $p$ is a prime and $m$ is a non-negative integer,
then there are at most two distinct integer squares among all $y_{k}>\max \left( 1,b^{2}/256 \right)$.
\end{proposition}

Part~(a) of this proposition has already been established. It is
Proposition~\genProp{} in \cite{V5} combined with Lemma~\gnLemma{} there
-- recall again that $N_{\alpha}$ there equals $-b$ here and $d$ there is
$a^{2}+b$ here. So we need only prove parts~(b) and (c).
The proof of this proposition is very similar to the proof of Proposition~\genProp{}
in \cite{V5}, so we refer to the proof there along with its intermediate steps
and inequalities there as much as possible. The proof is broken into steps,
following the steps and technique used to prove Proposition~\genProp{} in \cite{V5}.

\subsection{Proof Overview}
\label{subsect:overview}

We will assume the existence of a pair of distinct squares, $1<y_{k}<y_{\ell}$.
From this hypothetical pair, we obtain a very good approximation
to an algebraic number, $\omega_{k}^{1/4}$, by elements of an imaginary quadratic
field -- see equation~\eqref{eq:29} below. However, we can also construct
explicit approximations to $\omega_{k}^{1/4}$ using the classical (or Gauss)
hypergeometric functions, ${}_{2}F_{1}$.
These explicit approximations show that the approximation from our assumption
is ``too good''. In this way, in \cite{V5} we were able to obtain upper bounds
for $y_{\ell}$ in terms of $y_{k}$ (see equations~(\eqYEllUBa{}) and (\eqYEllUB{})
in \cite{V5}). We use these bounds, along with the lemmas in Section~\ref{sect:prelim},
to obtain a contradiction from our assumption, provided that $b \left( a^{2}+b \right)$
is sufficiently large. In Subsections~\ref{subsect:step-v} and \ref{subsect:step-vi},
we prove the proposition in these remaining cases. We also note here that \eqref{eq:yk-stepV-LB},
the main result in Subsection~\ref{subsect:step-v} may be of use for further work,
including work on particular equations.

In fact, we will initially assume there are three such distinct squares in the
sequence, as we will need this assumption to obtain the hypothetical pair of
squares with the required properties in Subsections~~\ref{subsect:step-i} and
~\ref{subsect:step-ii} below.

There is PARI/GP code available for checking all the numerical values in the proof.
See \verb!bPP-Prop31-sequence-proof-calcs.gp! in the \verb!bPP! directory of
\url{https://github.com/PV-314/hypgeom}.

\subsection{Choice of $k$ and $\ell$}
\label{subsect:k-and-ell}

To prove part~(a) of Proposition~\ref{prop:4.1}, we assume there are
two distinct squares, $y_{k}$ and $y_{\ell}$, in the sequence
with $y_{\ell}>y_{k}>1$ and we will obtain a contradiction.

To prove parts~(b) and (c) of Proposition~\ref{prop:4.1}, we assume there are
three distinct squares, $y_{k_{1}}$, $y_{k_{2}}$ and $y_{k_{3}}$, in the sequence
with $y_{k_{3}}>y_{k_{2}}>y_{k_{1}}>1$ and we will obtain a contradiction.
We let $k$ and $\ell$ be the two of the three indices, $k_{1}$, $k_{2}$ and $k_{3}$,
from the statement of Lemma~\ref{lem:gap}(b) (labelled $i$ and $j$, respectively there).

We put
$\omega_{k}=\left( x_{k} + N_{\varepsilon^{k}} \sqrt{-b} \right)
/\left( x_{k} - N_{\varepsilon^{k}} \sqrt{-b} \right)$. This is a corrected
version of the $\omega$ introduced in Section~4.1 of \cite{V5}. It is indexed by
$k$ to indicate its dependence on $k$ and we also correct, with no impact on our
argument, $N_{\varepsilon}$ there to $N_{\varepsilon^{k}}$. Similarly, the associated
$\varphi$ there is denoted by $\varphi_{k}$ here.

Let $\zeta_{4}$ be the $4$-th root of unity such that
\[
\left| \omega_{k}^{1/4} - \zeta_{4} \frac{x-y\sqrt{-\core(b)}}{x+y\sqrt{-\core(b)}} \right|
= \min_{0 \leq j \leq 3} \left| \omega_{k}^{1/4} - e^{2j\pi i/4} \frac{x-y\sqrt{-\core(b)}}{x+y\sqrt{-\core(b)}} \right|,
\]
where $x+y\sqrt{-\core(b)}=\left(  r_{k}-s_{k}\sqrt{-\core(b)} \right)
\left( r_{\ell}+s_{\ell}\sqrt{-\core(b)} \right)$ with $\left( r_{k}, s_{k} \right)$
and $\left( r_{\ell}, s_{\ell} \right)$ as in Proposition~\repProp{} of \cite{V5},
which are associated with $\left( x_{k}, y_{k} \right)$ and
$\left( x_{\ell}, y_{\ell} \right)$, respectively.
As in Subsection~\subSectPreq{} of \cite{V5}, we can take
$\zeta_{4} \in \bbQ \left( \sqrt{-\core(b)} \right)$. This is immediate for part~(a) of Proposition~\ref{prop:4.1}, where $b$ is a perfect square. For parts~(b) and (c),
this is true because in the notation of the proof of Lemma~\ref{lem:gap}(b), there are two $4$-th
roots of unity associated with $y_{k}$ and $y_{\ell}$ that satisfy
$\zeta_{4}^{(k)}=\pm \zeta_{4}^{(\ell)}$. Since
$\zeta_{4}^{(k)}=\pm \zeta_{4}^{(\ell)}$ and by Lemma~\miscLemma{}(b) of \cite{V5},
we have $\zeta_{4}=\pm 1$.

This is important for us here as $\zeta_{4} \left( x-y\sqrt{\core (b)} \right)
/ \left( x+y\sqrt{\core (b)} \right)$ must be in an imaginary quadratic field
in order to apply Lemma~\lemDio{} in \cite{V5} to obtain a lower bound for the
rightmost quantity in \eqref{eq:29} below.

\subsection{Further Prerequisites}
\label{subsect:preq}

As in the proof of Proposition~\genProp{} in \cite{V5}, we may assume that
$b \geq 2$. So for $b=p^{m}$, we may assume that $m>0$.

We also assume initially that
\begin{equation}
\label{eq:ab-LB}
a^{2}+b \geq 128.
\end{equation}

With the quantities defined in Subsection~\ref{subsect:k-and-ell}, equation~(\yEllLBEqnA) in \cite{V5} is
\begin{equation}
\label{eq:29}
\frac{2\sqrt{b}}{\sqrt{a^{2}+b} \, y_{\ell}}
= \left| \omega_{k} - \left( \frac{x-y\sqrt{-\core(b)}}{x+y\sqrt{-\core(b)}} \right)^{4} \right|
> 3.959 \left| \omega_{k}^{1/4} - \zeta_{4} \frac{x-y\sqrt{-\core(b)}}{x+y\sqrt{-\core(b)}} \right|.
\end{equation}

Since $y_{k} \geq 4$, we obtain
\begin{equation}
\label{eq:25}
x_{k}^{2}= \left( a^{2}+b \right) y_{k}^{2}-b
> \left( a^{2}+b \right) \left( y_{k}^{2}-1 \right)
\geq 0.9375 \left( a^{2}+b \right) y_{k}^{2},
\end{equation}
and so
\begin{equation}
\label{eq:26}
\sqrt{x_{k}^{2}+b}= \sqrt{\left( a^{2}+b \right) y_{k}^{2}}
< 1.04 \left| x_{k} \right|.
\end{equation}

We will use many of the same quantities here as in \cite{V5}. For each positive integer,
$r$, we denote by $p_{r}/q_{r}$,
the explicit approximations we mentioned in Subsection~\ref{subsect:overview},
defined in equation~(\eqApproxDef{}) of \cite{V5}, where $\omega$ there
is our $\omega_{k}$ above, $t'=-\core(b)$, $u_{1}=2x_{k}$, $u_{2}=2N_{\varepsilon^{k}} \sqrt{b/\core(b)}$
and $d'$ is as defined in equation~(\eqDPrimeDef{}) of \cite{V5}. Also, $k_{0}$,
$\ell_{0}$, $E$ and $Q$ are defined in Lemma~\lemDio{} of \cite{V5}.

From equation~(\eqELB) in \cite{V5}, we have
\begin{align}
\label{eq:E-LB1}
E
& > \frac{|g|\cN_{d',4}\left| \left| 2x_{k} \right| + 2\sqrt{x_{k}^{2}+b} \right|}{4e^{1.68}b}
> \frac{0.183|g|\cN_{d',4}\sqrt{a^{2}+b} \, y_{k}}{b} \\
& > 1.13, \nonumber
\end{align}
where $g$ is defined in equation~(\eqGDef{}) of \cite{V5} and $\cN_{d',4}$ is as
defined in equation~(\eqNdn{}) there with $n=4$. Hence $E>1$, as required for its
use with Lemma~\lemDio{} in \cite{V5}.

Similarly,
\begin{equation}
\label{eq:Q-LB1}
Q > \frac{2e^{1.68}\left( 1+\sqrt{0.9375} \right) \sqrt{\left( a^{2}+b \right)y_{k}^{2}}}{|g|\cN_{d',4}}
> 217
\end{equation}
and
\begin{equation}
\label{eq:Q-UB2}
Q<\frac{21.47\sqrt{a^{2}+b} \, y_{k}}{|g|\cN_{d',4}}.
\end{equation}

The numerical lower bounds in \eqref{eq:E-LB1} and \eqref{eq:Q-LB1} for $E$ and $Q$
are obtained using
Lemma~\EQLBLemma{} in \cite{V5}. That lemma requires $a^{2}+b \geq 105$. This
holds here due to \eqref{eq:ab-LB}.

We take $k_{0}=0.89$ (see equation~(\eqKZero{}) of \cite{V5}) and have
\begin{equation}
\label{eq:ell-UB}
\ell_{0}<0.458\sqrt{b}/ \left| x_{k} \right|.
\end{equation}

These four displayed inequalities, \eqref{eq:E-LB1}--\eqref{eq:ell-UB},
are equations~(\eqELB), (\eqQLB), (\eqQUB) and (\eqEllUB) of \cite{V5},
respectively.

Also from Lemma~\miscLemma{}(a) in \cite{V5}, we have $\left| \varphi_{k} \right|<0.6$,
so the condition $\left| \omega_{k}-1 \right|<1$ required in Lemma~\lemHypg{}\,
of \cite{V5} to apply the hypergeometric method is satisfied.

We let $B=x+y\sqrt{-\core(b)}=\left( r_{k}-s_{k}\sqrt{-\core(b)} \right) \left( r_{\ell}+s_{\ell}\sqrt{-\core(b)} \right)$
and $A=x-y\sqrt{-\core(b)}$. Recall equation~(\eqQAbs{}) of \cite{V5} (note that
these quantities are denoted by $q$ and $p$ respectively in \cite{V5}, but we have
already used $p$ in this paper),
\begin{equation}
\label{eq:qAbs}
\left| x+y\sqrt{-\core(b)} \right| = |B|
= \sqrt{f_{k}f_{\ell}} \left( y_{k}y_{\ell} \right)^{1/4}.
\end{equation}

As in \cite{V5}, we will break the proof into parts according to the value
of $r_{0}$, as defined in Lemma~\lemDio{} of \cite{V5}, the smallest positive
integer such that
$\left( Q-1/E  \right)\ell_{0}|B|/\left( Q-1 \right)<cE^{r_{0}}$.
For the same reasons as in Subsection~\subsectStepiii{} of \cite{V5},
we will take $c=0.75$ here and throughout the remainder of Section~\ref{sect:prop-11}.

\subsection{$r_{0}=1$ and $\zeta_{4}A/B \neq p_{1}/q_{1}$ for all $4$-th roots of unity, $\zeta_{4}$}
\label{subsect:step-i}

We start by determining an upper bound for $y_{\ell}$ for all $r_{0} \geq 1$ when
$\zeta_{4}A/B \neq p_{r_{0}}/q_{r_{0}}$, since we will also need such a result
in Subsection~\ref{subsect:step-iii}.

Recall equation~(\eqYEllUBa) from \cite{V5}:
\[
\left( bf_{k}f_{\ell} \right)^{2}\left( \frac{0.45}{1-c} \right)^{4}
\left( \frac{21.47}{|g|\cN_{d',4}} \right)^{4r_{0}}
\left( a^{2}+b \right)^{2r_{0}-2} y_{k}^{4r_{0}+1}
>y_{\ell}^{3}.
\]

With $c=0.75$ as above, if this inequality holds, then so does
\begin{equation}
\label{eq:y2UB-step1}
10.5\left( bf_{k}f_{\ell} \right)^{2}
\left( \frac{21.47}{|g|\cN_{d',4}} \right)^{4r_{0}}
\left( a^{2}+b \right)^{2r_{0}-2} y_{k}^{4r_{0}+1}
>y_{\ell}^{3}.
\end{equation}

Specialising to the case when $r_{0}=1$ and using $|g|\cN_{d',4} \geq 2$ from
Lemma~\gnLemma{} of \cite{V5}, we must have
\begin{equation}
\label{eq:y3UB-step1a}
y_{\ell}^{3} < 140000 \left( bf_{k}f_{\ell} \right)^{2}y_{k}^{5}.
\end{equation}

This inequality is independent of the choice of $y_{k}$ and $y_{\ell}$,
provided $y_{\ell} \neq y_{k}$ and $y_{\ell}, y_{k}>1$.

Combining Lemma~\ref{lem:gap}(b) with \eqref{eq:y3UB-step1a}, we have
\begin{equation}
\label{eq:y3UB-step1b}
140000 \left( bf_{k}f_{\ell} \right)^{2}y_{k}^{5}>y_{\ell}^{3}
>15^{3} \left( \frac{a^{2}+b}{b} \right)^{6} y_{k}^{9}.
\end{equation}

We have $f_{k}, f_{\ell} \leq 2$
when $b$ is even. Applying this and \eqref{eq:ab-LB} to \eqref{eq:y3UB-step1b}, we obtain
\[
\frac{b^{2}}{285} >
\frac{140000^{1/4} \cdot 16^{1/4}}{15^{3/4}} \frac{b^{2}}{\left( a^{2}+b \right)^{3/2}}
>y_{k}.
\]

But from Lemma~\ref{lem:Y-LB}(a) and (b), this is not possible.

Similarly, when $b$ is odd, we have $f_{k}, f_{\ell} \leq 4$ and obtain
\[
y_{k} < \frac{b^{2}}{142},
\]
which is again not possible by Lemma~\ref{lem:Y-LB}.

Hence we cannot have three distinct squares all greater than $1$ in this case.

\subsection{$r_{0}=1$ and $\zeta_{4}A/B = p_{1}/q_{1}$ for some $4$-th root of
unity, $\zeta_{4}$}
\label{subsect:step-ii}

In Subsection~\rEqualSubsection{} of \cite{V5} (see equation~(\eqYEllUB) there),
we showed that
\begin{equation}
\label{eq:y2UB-step2}
1.734r_{0}^{1/2} \left(4\frac{a^{2}+b}{b} \right)^{r_{0}}y_{k}^{2r_{0}+1}
>y_{\ell}.
\end{equation}

We again specialise to the case of $r_{0}=1$.

Our gap principle in Lemma~\ref{lem:gap} and \eqref{eq:y2UB-step2} imply that
\begin{equation}
\label{eq:yUB-step2}
6.94 \frac{a^{2}+b}{b} y_{k}^{3}>y_{\ell}
>15 \frac{\left( a^{2}+b \right)^{2}}{b^{2}} y_{k}^{3}.
\end{equation}
But this inequality is impossible.

\subsection{$r_{0}>1$, $\zeta_{4}A/B \neq p_{r_{0}}/q_{r_{0}}$ for all $4$-th
roots of unity, $\zeta_{4}$, and large $a^{2}+b$}
\label{subsect:step-iii}

Here we establish a stronger gap principle for $y_{k}$ and $y_{\ell}$ than the
one in Lemma~\ref{lem:gap}. We then use this with the upper bound for $y_{\ell}$
in \eqref{eq:y2UB-step1} to obtain a contradiction.

As in equation~(\eqnYKUBa) of \cite{V5}, we have
\begin{equation}
\label{eq:y1UB-step3-b}
\left( f_{k}f_{\ell} \right)^{8}
      \frac{0.00078b^{2} \left( |g|\cN_{d',4} \right)^{4}}{\left( a^{2}+b \right)^{4}}
      \left( \frac{1.22 \cdot 10^{7}}
                  {\left( |g|\cN_{d',4} \right)^{8}} \right)^{2r_{0}-1}
> \left( \frac{y_{k}^{4}\left( a^{2}+b \right)^{2}}{b^{6}} \right)^{2r_{0}-1}.
\end{equation}

Suppose $b$ is even, we have $f_{k}, f_{\ell} \leq 2$ and $|g|\cN_{d',4} \geq 2\sqrt{2}$
(by Lemma~\gnLemma{} in \cite{V5}).
For $r_{0}>1$, $x^{12-16r_{0}}$ is a decreasing
function of $x$ for $x>1$. So \eqref{eq:y1UB-step3-b} can only hold if
\[
\frac{3320b^{2}}{\left( a^{2}+b \right)^{4}}
> \left( \frac{y_{k}^{4}\left( a^{2}+b \right)^{2}}{2980b^{6}} \right)^{2r_{0}-1}
> \left( \frac{b^{2}\left( a^{2}+b \right)^{2}}{1.28 \cdot 10^{13}} \right)^{2r_{0}-1}.
\]
holds, where the last inequality comes from applying $y_{k} \geq b^{2}/256$
in Lemma~\ref{lem:Y-LB}.

By \eqref{eq:ab-LB}, the left-hand side is less than $1$. But if
\begin{equation}
\label{eq:step-iii-bEvenLB}
b \left( a^{2}+b \right) > \sqrt{1.28 \cdot 10^{13}}=3.577 \ldots \cdot 10^{6},
\end{equation}
then the right-hand side is greater than $1$. Hence for $b$ even, if
\eqref{eq:step-iii-bEvenLB} holds, then this case (namely,
$r_{0}>1$ and $\zeta_{4}A/B \neq p_{r_{0}}/q_{r_{0}}$ for all $4$-th
roots of unity, $\zeta_{4}$) is not possible.

Similarly, if $b$ is an odd prime power, we apply $f_{k},f_{\ell} \leq 4$,
$|g|\cN_{d',4} \geq 2$ (by Lemma~\gnLemma{} in \cite{V5})
and $y_{k} \geq b^{2}/16$ to \eqref{eq:y1UB-step3-b}. Thus
\begin{equation}
\label{eq:ykUB-step3-bOdd}
\frac{5.43 \cdot 10^{7}b^{2}}{\left( a^{2}+b \right)^{4}}
> \left( \frac{y_{k}^{4}\left( a^{2}+b \right)^{2}}{47700b^{6}} \right)^{2r_{0}-1}
> \left( \frac{b^{2}\left( a^{2}+b \right)^{2}}{3.13 \cdot 10^{9}} \right)^{2r_{0}-1}.
\end{equation}

If
\begin{equation}
\label{eq:step-iii-bOddLB1}
b\left( a^{2}+b \right) > 5.43 \cdot 10^{7},
\end{equation}
then $\left( a^{2}+b \right)^{4}/b^{2}> 5.43 \cdot 10^{7}$ and the left-hand side
of \eqref{eq:ykUB-step3-bOdd} is less than $1$. But if
\begin{equation}
\label{eq:step-iii-bOddLB2}
b \left( a^{2}+b \right) > \sqrt{3.13 \cdot 10^{9}}=55946.4\ldots,
\end{equation}
then the right-hand side of \eqref{eq:ykUB-step3-bOdd} is greater than $1$.
Hence for $b$ odd, if \eqref{eq:step-iii-bOddLB1} holds, then this case is not
possible.

We will treat the outstanding values of $a$ and $b$ in Subsection~\ref{subsect:step-vi}
below.

\subsection{$r_{0}>1$, $\zeta_{4}A/B = p_{r_{0}}/q_{r_{0}}$ for some $4$-th
root of unity, $\zeta_{4}$, and large $a^{2}+b$}
\label{subsect:step-iv}

As in equation~(\eqnYKUBb) of \cite{V5}, we have
\begin{equation}
\label{eq:step-iv}
1 > \frac{0.7405}{f_{k}^{2}f_{\ell}^{2}}
  \left( \frac{0.0002344 |g|^{4}\cN_{d',4}^{4}}{b^{3}} \left( a^{2}+b \right)y_{k}^{2} \right)^{r_{0}-1}.
\end{equation}

If $b$ is even, then we have $f_{k}, f_{\ell} \leq 2$,
$|g|\cN_{d',4} \geq 2\sqrt{2}$ and $y_{k}>b^{2}/256$ by Lemma~\ref{lem:Y-LB}(a).
So if \eqref{eq:step-iv} holds, then so does
\begin{align}
\label{eq:ykUB-step4-bEven}
1 & > 0.046
\left( \frac{0.015y_{k}^{2}\left( a^{2}+b \right)}{b^{3}} \right)^{r_{0}-1} \\
& >\left( \frac{0.00068y_{k}^{2}\left( a^{2}+b \right)}{b^{3}} \right)^{r_{0}-1}
>\left( 10^{-8}b\left( a^{2}+b \right) \right)^{r_{0}-1}. \nonumber
\end{align}

But this is not possible if
\begin{equation}
\label{eq:step-iv-bEven-LB}
b \left( a^{2}+b \right) \geq 10^{8}.
\end{equation}

Hence for $b$ even, if \eqref{eq:step-iv-bEven-LB} holds, then this case
(namely, $r_{0}>1$ and $\zeta_{4}A/B = p_{r_{0}}/q_{r_{0}}$ for some $4$-th root
of unity, $\zeta_{4}$) is not possible.

If $b$ is odd, then we have $f_{k}, f_{\ell} \leq 4$,
$|g|\cN_{d',4} \geq 2$ and $y_{k}>b^{2}/16$ by Lemma~\ref{lem:Y-LB}(a).
So if \eqref{eq:step-iv} holds, then so does
\begin{align}
\label{eq:step-iv-ykUB-bOdd}
1 & > 0.0028
\left( \frac{0.0037y_{k}^{2}\left( a^{2}+b \right)}{b^{3}} \right)^{r_{0}-1} \\
&> \left( \frac{0.00001y_{k}^{2}\left( a^{2}+b \right)}{b^{3}} \right)^{r_{0}-1}
> \left( 3.9 \cdot 10^{-8}b\left( a^{2}+b \right) \right)^{r_{0}-1}. \nonumber
\end{align}

But this is not possible if
\begin{equation}
\label{eq:step-iv-bOdd-LB}
b \left( a^{2}+b \right) \geq 26,000,000.
\end{equation}

Therefore, for $b$ odd, if \eqref{eq:step-iv-bOdd-LB} holds, then this case is not
possible.

We will treat the outstanding values of $a$ and $b$ in Subsection~\ref{subsect:step-vi}
below.

\subsection{$r_{0}>1$ and large $y_{k}$}
\label{subsect:step-v}

We will show that for $r_{0}>1$ we cannot have
\begin{equation}
\label{eq:yk-stepV-LB}
y_{k} > 1800b^{3/2}.
\end{equation}

From \eqref{eq:yk-stepV-LB} and \eqref{eq:25},
\[
x_{k}^{2}>0.9375 \left( a^{2}+b \right) y_{k}^{2}
> 0.9375 \left( a^{2}+b \right) \left( 1800b^{3/2} \right)^{2}
> 3 \cdot 10^{6} b^{3}.
\]

For such $x_{k}$, from \eqref{eq:E-LB1} and $|g|\cN_{d',4} \geq 2$, we have
\[
E^{2}> \left( \frac{2 \left| \left| 2x_{k} \right| + 2\sqrt{x_{k}^{2}+b} \right|}{4e^{1.68}b} \right)^{2}
>\frac{64x_{k}^{2}}{16e^{3.36}b^{2}}
>\frac{4 \cdot 3 \cdot 10^{6} b^{3}}{e^{3.36}b^{2}}>410,000b.
\]

From \eqref{eq:E-LB1}, \eqref{eq:Q-UB2} and $|g|\cN_{d',4} \geq 2$, we also have
\[
\frac{Q}{E}<\left( 10.74\sqrt{a^{2}+b} \, y_{k} \right)
\frac{b}{0.366\sqrt{a^{2}+b} \, y_{k}}
< 29.4b.
\]

Thus $E^{2}>Q/E$ and so $E^{3}>Q$.
Therefore, from $E>1$, the definition of $r_{0}$ (immediately before Subsection~\ref{subsect:step-i}),
\eqref{eq:Q-LB1}, \eqref{eq:ell-UB}, $c=0.75$ and \eqref{eq:qAbs},
\begin{align}
\label{eq:step5QUB}
Q^{r_{0}-2} &\leq E^{3(r_{0}-2)} = E^{-3}E^{3(r_{0}-1)}
\leq E^{-3} \left( \frac{(Q-1/E)\ell_{0}|B|}{c(Q-1)} \right)^{3} \\
&<E^{-3} \left( \frac{217}{216c} \frac{0.458\sqrt{b}}{\left| x_{k} \right|} |B| \right)^{3}
<E^{-3} \left( 0.614\frac{\sqrt{bf_{k}f_{\ell}}}{\left| x_{k} \right|} \left( y_{k}y_{\ell} \right)^{1/4} \right)^{3}. \nonumber
\end{align}
Here we use the numeric lower bound for $Q$ in \eqref{eq:Q-LB1} to bound $(Q-1/E)/(Q-1)$
from above. As noted after
\eqref{eq:Q-LB1}, this requires \eqref{eq:ab-LB}.

From \eqref{eq:29}, Lemma~\lemDio{} of \cite{V5}, \eqref{eq:qAbs}, $c=0.75$ and
$k_{0}=0.89$, we have
\begin{align*}
\frac{2\sqrt{b}}{\sqrt{a^{2}+b} \, y_{\ell}}
& > 3.959 \left| \omega_{k}^{1/4} - \zeta_{4} \frac{x-y\sqrt{-\core(b)}}{x+y\sqrt{-\core(b)}} \right|
> \frac{3.959(1-c)}{k_{0}Q^{r_{0}+1} \left| x+y\sqrt{-\core(b)} \right|} \\
& > \frac{1.11}{Q^{r_{0}+1} \sqrt{f_{k}f_{\ell}} \left( y_{k}y_{\ell} \right)^{1/4}},
\end{align*}
where $x$, $y$ and $\zeta_{4}$ are defined before \eqref{eq:29}.

Applying \eqref{eq:step5QUB} to this inequality, we obtain
\[
\frac{2\sqrt{b}}{\sqrt{a^{2}+b} \, y_{\ell}}
>\frac{1.11 \left| x_{k} \right|^{3}}{Q^{3}E^{-3} \left( 0.614\sqrt{bf_{k}f_{\ell}} \left( y_{k}y_{\ell} \right)^{1/4} \right)^{3}\sqrt{f_{k}f_{\ell}}\left( y_{k}y_{\ell} \right)^{1/4}}.
\]

We saw above that $Q/E<29.4b$, so
\[
\frac{2\sqrt{b}}{\sqrt{a^{2}+b}}
>\frac{1.11 \left| x_{k} \right|^{3}}{29.4^{3}b^{3} \left( 0.614\sqrt{bf_{k}f_{\ell}} y_{k}^{1/4} \right)^{3} \sqrt{f_{k}f_{\ell}} y_{k}^{1/4}}
>\frac{0.00018 \left| x_{k} \right|^{3}}{b^{9/2}f_{k}^{2}f_{\ell}^{2} y_{k}}.
\]

Squaring the far-left and far-right sides of this inequality, simplifying the
resulting inequality, then combining the result with \eqref{eq:25}, we find that
\[
124 \cdot 10^{6} \left( f_{k}f_{\ell} \right)^{4} b^{10} y_{k}^{2}> \left( a^{2}+b \right) x_{k}^{6}
> \left( a^{2}+b \right) \left( 0.9375 \left( a^{2}+b \right) y_{k}^{2} \right)^{3}.
\]

Applying $f_{k}f_{\ell} \leq 16$ to this inequality and simplifying, we obtain
\[
1800^{4}b^{6}>10^{13} \frac{b^{10}}{\left( a^{2}+b \right)^{4}}>y_{k}^{4},
\]

But this contradicts our lower bound for $y_{k}$ in \eqref{eq:yk-stepV-LB}.

\subsection{Small $a^{2}+b$ and small $y_{k}$}
\label{subsect:step-vi}

To complete the proof of Proposition~\ref{prop:4.1}, we need to remove the
assumption on $a^{2}+b$ in \eqref{eq:ab-LB} and those on $b \left( a^{2}+b \right)$
in Subsections~\ref{subsect:step-iii} and \ref{subsect:step-iv}
(see equations~\eqref{eq:step-iii-bEvenLB}, \eqref{eq:step-iii-bOddLB1},
\eqref{eq:step-iv-bEven-LB} and \eqref{eq:step-iv-bOdd-LB}).

\subsubsection{Equations~\eqref{eq:step-iii-bEvenLB}, \eqref{eq:step-iii-bOddLB1},
\eqref{eq:step-iv-bEven-LB} and \eqref{eq:step-iv-bOdd-LB}}

To handle the cases arising from these equations, we proceeded as follows.

For each positive integer $d$ satisfying $128 \leq d \leq 10^{8}$,
we checked each pair of positive integers
$(a,b)$ such that $a^{2}+b=d$,
$b\left( a^{2}+b \right) \leq 10^{8}$,
$b=p^{m}$, $2p^{m}$, $4p^{m}$, $8p^{m}$ or $16p^{m}$,
where $p$ is a prime and $m$ a non-negative integer.
We searched for squares in the sequence of $y_{k}$'s
satisfying $\max \left( 1, b^{2}/16 \right) < y_{k} < 1800b^{3/2}$ or
$\max \left( 1, b^{2}/256 \right) < y_{k} < 1800b^{3/2}$, as appropriate. The lower
bounds coming from Proposition~\ref{prop:4.1}(b) and (c) and the upper
bound from \eqref{eq:yk-stepV-LB}.
We used PARI/GP \cite{Pari} to perform this search. It took $13,978$
seconds on a laptop with an Intel Core i7-13700H@2.40GHz processor and 32~GB of RAM.

This PARI/GP code is in \verb!bPP-end-proof-check.gp! in the \verb!bPP! directory
of \url{https://github.com/PV-314/hypgeom}. See in particular, the function,
\verb!search_all_d()!.

This left 531 sequences. The largest value of $d$ among them was $853,864$. For
each of these sequences, there are associated values of $b$ and $d=a^{2}+b$, and
a diophantine equation, $X^{2}-dY^{4}=-b$. All squares in each sequence will be
among the solutions in positive integers of the associated diophantine equation.
We solved these diophantine equations using Magma (version V2.28-2)
\cite{Magma} and its \verb+IntegralQuarticPoints()+ function.
The total time taken by Magma's online calculator was 403 seconds. No counterexamples
to Proposition~\ref{prop:4.1} were found in this way for these sequences.

Here is further information about the results of these calculations.

All of the equations, $X^{2}-dY^{4}=-b$, have at most $4$ solutions in positive integers.
Only for
$(b, d) = (136, 152)$ and $(152, 168)$ are there 4 such solutions.
Only for
$(b,d)=(25,986)$, $(81,442)$, $(100,3944)$, $(324,1768)$ and $(625,24650)$ is
$b$ a square and did the equation have $3$ such solutions.

For $(b,d)=(81,442)$, i.e., the equation $X^{2}-442Y^{4}=-81$, note that
$\left( 189+3^{2}\sqrt{442} \right)/\left( 19+\sqrt{442} \right)
=\left( 43+2\sqrt{442} \right)/9$, so the solutions $(189,3)$ and $(19,1)$
arise from different sequences. We did the same check for the solutions of the
other equations. No exceptions to Proposition~\ref{prop:4.1} were found.

\subsubsection{Equation~\eqref{eq:ab-LB}}

Here we proceeded similarly, but simply checked all pairs $(a,b)$ satisfying
the conditions in Proposition~\ref{prop:4.1} and $a^{2}+b<128$.
There were 645 such pairs. The associated equations were solved using Magma,
as above. The total time taken by Magma's online calculator was 164 seconds.
Once again, no counterexamples to Proposition 3.1 were found in this way
for these sequences.

Here is further information about the results of these calculations.

All of the equations had at most $4$ solutions in positive integers.
Only for
$(b,d)=(16,17)$ and $(64,68)$ was $b$ a square and the equation
had $3$ such solutions.
Only for
$(b,d)=(34,38)$ and $(38,42)$
were there $4$ such solutions.
The same checks as illustrated above for $(b,d)=(81,442)$ showed that these do
not give rise to exceptions to Proposition~\ref{prop:4.1}.

This completes the
proof of Proposition~\ref{prop:4.1}.

\section{Proofs of Theorems and Corollaries}
\label{sect:proofs}

\subsection{Proof of Theorem~\ref{thm:1.2-seq}}

If $b$ is an even square, then by Proposition~\ref{prop:4.1}(a), there is at most
one square, $y_{k}$, satisfying $y_{k}>4.75b^{3/2}/\sqrt{a^{2}+b}$. Applying
Lemma~\ref{lem:Y-LB}, we see
that if $b^{2}/256>4.75b^{3/2}/\sqrt{a^{2}+b}$, then the theorem holds. This will
be the case if $b \left( a^{2}+b \right)>1.5 \cdot 10^{6}$.

We proceed as in Subsection~\ref{subsect:step-vi}, searching for pairs $(a,b)$
such that $b=2^{m}$, $4p^{m}$ or $16p^{m}$ is a square and $b \left( a^{2}+b \right) \leq 1.5 \cdot 10^{6}$.
For each such pair, $(a,b)$, we let $(t,u)$ be the smallest solution in positive
integers of $X^{2}-\left( a^{2}+b \right) Y^{2}=\pm 4$ and calculate elements
of the sequence $\left( y_{k} \right)_{k=-\infty}^{\infty}$ satisfying
$\max \left( 1, b^{2}/256 \right) < y_{k} \leq 4.75b^{3/2}/\sqrt{a^{2}+b}$.
Five such sequences were found, for all but one of them $a^{2}+b \leq 128$, so
they have already been tested by the work in Subsection~\ref{subsect:step-vi}.
The remaining sequence, formed from $a=10$, $b=10^{2}$, $t=14$ and $u=1$, is checked
by using Magma's \verb+IntegralQuarticPoints()+ function to find all integer
solutions of \eqref{eq:2}. This calculation completed the proof of the theorem
when $b$ is an even square.

We proceed in the same way when $b$ is an odd square. We find two sequences in
this case, but for both of them $a^{2}+b \leq 128$, so they have already been treated.

When $b$ is not a square, Theorem~\ref{thm:1.2-seq} follows immediately from
Proposition~\ref{prop:4.1}(b) and (c).

\subsection{Proof of Corollary~\ref{cor:1.2-eqn1}}

The equation~\eqref{eq:2} has the solution $(X,Y)=(a,1)$. So we will apply
Lemma~\ref{lem:pell-rep} with $d=a^{2}+b$, $e=-b$, $x_{0}=a$ and $y_{0}=1$.
The assumption that
$\gcd \left( a^{2}, b \right)$ is squarefree implies that $\gcd \left( d_{2}, e \right)=1$
in Lemma~\ref{lem:pell-rep} holds. We also assume that $d$ is not a square.

Writing $d=d_{1}d_{2}^{2}=a^{2}+b$, if $d_{1} \equiv 1 \pmod{8}$, then
either $d \equiv 1 \pmod{8}$ or $d \equiv 0 \pmod{4}$. The case when
$d \equiv 1 \pmod{8}$ is clearly covered by the hypotheses of the corollary.
When $d \equiv 0 \pmod{4}$, then $2 | d_{2}$. Hence $b=-e$ must be odd by the
hypothesis in Lemma~\ref{lem:pell-rep} that $\gcd \left( d_{2}, e \right)=1$.
This case is not excluded by the hypotheses of Lemma~\ref{lem:pell-rep}. Hence
Lemma~\ref{lem:pell-rep} is applicable under the conditions in Corollary~\ref{cor:1.2-eqn1}.

So $\eta$ in Lemma~\ref{lem:pell-rep} can be written as
$\left( t'+u'\sqrt{d} \right)/2$, a unit of norm $1$.
As $e<0$, we also have $\eta>0$ and that $\eta$ is a power of the fundamental
unit of $\cO_{\bbQ\left( \sqrt{a^{2}+b} \right)}$ from Lemma~\ref{lem:pell-rep}.
Since $X^{2}-dY^{2}=-4$ is solvable, it follows that the fundamental unit of
$\cO_{\bbQ\left( \sqrt{a^{2}+b} \right)}$
has norm $-1$. Thus $\eta$ is an even power of a unit of the form $\varepsilon=\left( t+u\sqrt{d} \right)/2$
with norm $-1$. Hence by Lemma~\ref{lem:pell-rep}, all other solutions of \eqref{eq:2}
must come from squares in the single sequence of $y_{k}$'s defined in \eqref{eq:yk-defn}
for this $\varepsilon$. Corollary~\ref{cor:1.2-eqn1} now follows by applying
Theorem~\ref{thm:1.2-seq} to this sequence.

\subsection{Proof of Corollary~\ref{cor:1.2-eqn2}}

For $m=1$, $2$ and $3$, if $\gcd(X,Y) \neq 1$, then we must have
$\gcd(X,Y)=p$ by \eqref{eq:2}. This immediately implies that $m=1$ is not possible,
since $p^{2}$ divides the left-hand side of \eqref{eq:2}.
When $m=3$, it implies that $p^{3}$, and hence $p^{4}$, divides $X^{2}$.
But this means that $p^{4}$ divides the left-hand side of \eqref{eq:2}, which is
also not possible. So we only have the at most three coprime solutions from
Corollary~\ref{cor:1.2-eqn1}(b) in these two cases.

For $m=2$, we can remove a common factor of $p^{2}$ from each term in
\eqref{eq:2}, obtaining $X_{1}^{2}- \left( a^{2}+b \right)p^{2}Y_{1}^{4}=-1$,
where $X=pX_{1}$ and $Y=pY_{1}$. Since $\left( a^{2}+b \right)p^{2}=\left( a^{2}+p^{2} \right)p^{2}
\geq 20$, we can now appeal to Theorem~D of \cite{Chen1} (using $\left( a^{2}+b \right)p^{2}$
as $d$ there) to show there is at most one such solution. Combining this with
the at most two coprime solutions from Corollary~\ref{cor:1.2-eqn1}(a),
we obtain the corollary in this case.

For $m=4$, we have $p^{4}|X^{2}$ by the same argument as for $m=3$. Now we can
remove a common factor of $p^{4}$ from each term in \eqref{eq:2}, obtaining
$X_{1}^{2}- \left( a^{2}+b \right)Y_{1}^{4}=-1$, where $X=p^{2}X_{1}$ and
$Y=pY_{1}$. We now conclude as for $m=2$.

\section{Examples}
\label{sect:egs}

In this section, we give examples showing that our results are best possible.

The following example will be used several times in this section.

\begin{example}
\label{eg:1}
Let $b>5$ satisfy $b \equiv 1 \pmod{4}$ and not divisible by $5$ and
$a=\left( b-5 \right)/4$.
With $t=2a+4$ and $u=2$, $\varepsilon=\left( t+u\sqrt{a^{2}+b} \right)/2$
is a unit in $\cO_{\bbQ \left( \sqrt{a^{2}+b} \right)}$ and
we have $y_{1}=(b+1)^{2}/4$.
\end{example}

PARI/GP code for finding the examples below is in the files\\
\verb!bPP-eqn-solution-search.gp!
and \verb!bPP-square-search.gp! in the \verb!bPP! directory of \url{https://github.com/PV-314/hypgeom}.

\subsection{Examples for Theorem~\ref{thm:1.2-seq} and Corollary~\ref{cor:1.2-eqn1} with $b=p^{m}$}
\label{subsect:thm12-pm-egs}

Let $b \geq 9$ be the square of an odd prime power not divisible by $5$. Then
Example~\ref{eg:1} shows that Theorem~\ref{thm:1.2-seq} and Corollary~\ref{cor:1.2-eqn1}
are best possible if $b$ is a square.

When $b$ is an arbitrary prime power, Theorem~\ref{thm:1.2-seq} and
Corollary~\ref{cor:1.2-eqn1} are also best possible.
We found three examples where there are three distinct squares in the sequences,\\
-- $(a,b,t,u)=\left( 1,7,2,1 \right)$, so that $y_{0}=1$, $y_{1}=2^{2}$ and
$y_{-3}=8^{2}$;\\
-- $(a,b,t,u)=\left(1,11,4,1 \right)$, so that $y_{0}=1$, $y_{1}=3^{2}$ and
$y_{-3}=31^{2}$;\\
-- $(a,b,t,u)=\left(36,29,364,10 \right)$, so that
$y_{0}=1$, $y_{1}=363^{2}$ and $y_{-1}=27^{2}$.

It is conceivable that there are only two distinct integer squares for all other $a$ and
$b$ satisfying the conditions of Theorem~\ref{thm:1.2-seq} with $b=p^{m}$ not
a square. If so, this would be best possible by the family of examples
above.

\subsection{Examples for Theorem~\ref{thm:1.2-seq} with $b=2p^{m}$}
\label{subsect:thm12-2pm-egs}

Theorem~\ref{thm:1.2-seq} is best possible when $b=2p^{m}$ is a square too. We
found one example: $(a,b,t,u)=\left( 1,2 \cdot 2^{3}, 8,2 \right)$ and
$y_{-1}=5^{2}$.

When $b=2p^{m}$ is not square, our results
do not appear to be best possible. We found a family of sequences such that
$y_{-1}$ is also a square. Let $a$ be an even integer, put $b=3a^{2}-8a+6$,
$t=8a^{2}-16a+10$ and $u=4a-4$. Then $\varepsilon$ has norm $1$ and
$y_{-1}=\left( 4a^{2}-10a+7 \right)^{2}$. There are conjecturally infinitely
many integers $a$ such that $b$ is twice a prime.

Since $\varepsilon$ has norm $+1$, rather than $-1$, this family of
examples does not provide examples for Corollary~\ref{cor:1.2-eqn1}.

We found three examples with $b=2p^{m}$ is not square, $\gcd(a,b)=1$ and two
integer squares:\\
-- $(a,b,t,u)=\left( 3, 2 \cdot 2^{2}, 8,2 \right)$ and $y_{-1}=3^{2}$;\\
-- $(a,b,t,u)=\left( 3, 2 \cdot 2^{4}, 64,10 \right)$ and
$y_{-1}=33^{2}$; and\\
-- $(a,b,t,u)=\left(45, 2 \cdot 2^{6}, 464,10 \right)$ and $y_{-1}=57^{2}$.

It is possible that there are no integer squares $y_{k} \neq 1$ for all other
sequences formed from relatively prime $a$ and $b$ satisfying the conditions of
Theorem~\ref{thm:1.2-seq} with $b=2p^{m}$ not a square.

\subsection{Examples for Theorem~\ref{thm:1.2-seq} and Corollary~\ref{cor:1.2-eqn1} with $b=4p^{m}$}
\label{subsect:thm12-4pm-egs}

Theorem~\ref{thm:1.2-seq} is best possible when $b=4p^{m}$ is a square. We have
an infinite family of examples to
demonstrate this. Let $a$ be any positive integer and put $b=8a+20$, so
$d=a^{2}+8a+20$ and we can use $t=a+4$ and $u=1$. Then $y_{1}=(a+3)^{2}$. If we
also suppose that $b=4p^{m}$ for some prime $p$ and $m$ a positive even integer,
then we have the desired examples.

Note that $t^{2}-du^{2}=-4$, but $x_{1}=a^{3}+10a^{2}+35a+40$,
so when $a$ is odd, this gives rise to a solution of \eqref{eq:2} where both $X$
and $Y$ are even. If $a$ is even, then $4|\gcd \left( a^{2}, b \right)$, so the
condition in Corollary~\ref{cor:1.2-eqn1} is not satisfied. Hence this does not give
a family of examples for Corollary~\ref{cor:1.2-eqn1}.

We have two examples showing that Corollary~\ref{cor:1.2-eqn1} is best possible
when $b=4p^{m}$ is a square:\\
--for $(a,b,d)= \left( 3,484=4\cdot 11^{2}, 493 \right)$, we have the solutions
$(X,Y)=(3, 1)$ and $(118323, 73)$,\\
--for $(a,b,d)= \left( 3, 676=4 \cdot 13^{2}, 685 \right)$, we have the solutions
$(X,Y)=(3, 1)$ and $(6674643, 505)$.

When $b=4p^{m}$ is not square, Theorem~\ref{thm:1.2-seq} is best possible.
We found one example:\\
-- $(a,b,t,u)=\left( 72,116, 364,5  \right)$, so that $y_{1}=363^{2}$ and
$y_{-1}=27^{2}$.

\subsection{Examples for Theorem~\ref{thm:1.2-seq} with $b=8p^{m}$}
\label{subsect:thm12-8pm-egs}

When $p$ is an odd prime, $b$ is never a square. For such $b$, we have not found
any sequences with three distinct squares, although one can construct infinite
families of such sequences with two distinct squares.

When $p=2$, see Subsection~\ref{subsect:thm12-2pm-egs} for examples.

\subsection{Examples for Theorem~\ref{thm:1.2-seq} with $b=16p^{m}$}
\label{subsect:thm12-16pm-egs}

When $b=16p^{m}$ is a square, Theorem~\ref{thm:1.2-seq} is best possible.

We define a sequence $\left( u_{n} \right)_{n \geq 0}$ by $u_{0}=6$,
$u_{1}=88$ and $u_{n+2}=14u_{n+1}-u_{n}+4$ for $n \geq 2$. We also define
$\left( a_{n} \right)_{n \geq 0}$ by $a_{n}=4u_{n}-4$
and
$\left( t_{n} \right)_{n \geq 0}$ by $t_{n}=8u_{n}^{2}+2$.
Lastly, define $\left( b_{1,n} \right)_{n \geq 0}$ by $b_{1,0}=11$,
$b_{1,1}=153$ and $b_{1,n+2}=14b_{1,n+1}-b_{1,n}'$ for $n \geq 2$ and put
$b_{n}=16b_{1,n}^{2}$. The $b_{1,n}$'s and $u_{n}$'s are related via
$\left( 3u_{n}+1 \right)^{2}-3b_{1,n}^{2}=-2$. So
$d_{n}=a_{n}^{2}+b_{n}=64u_{n}^{2}+32$ and the norm of
$\left( t_{n}+u_{n} \sqrt{d_{n}} \right)/2$ is $+1$. Also
$y_{n,-1}=\left( 1+2u_{n}+4u_{n}^{2} \right)^{2}$.

Conjecturally, there should be infinitely many primes among the $b_{1,n}$'s.
For $n \leq 2500$, we found that $b_{1,n}$ was prime for
$n=3,9,37,69$, $79$, $298$ and $483$. 

We found no examples with $b=16p^{m}$ not a square showing that Theorem~\ref{thm:1.2-seq}
is best possible in this case, although we did find infinitely many such sequences
with two distinct squares.

\subsection{Examples for Corollary~\ref{cor:1.2-eqn2}}
\label{subsect:cor12-egs}

For $b$ prime, we have the three examples given above in Subsection~\ref{subsect:thm12-pm-egs}
with three solutions.

We found one example with $b$ the square of a prime and three solutions, namely
$(a,b)=(31,25)$ with the solutions $(X,Y)=(31, 1)$, $(785, 5)$ and $(3076289, 313)$.

No examples with three solutions were found for $b$ a cube.

For $b=p^{4}$, where $p>5$, we can use Example~\ref{eg:1}. In addition to the
two solutions given above, we also have the solution
$(X,Y)=\left( \left( p^{6}+3p^{2} \right)/4, p \right)$.

\subsection{Examples for Lemma~\ref{lem:Y-LB}}
\label{subsect:YLB-egs}

PARI/GP code for finding the examples below is in the file \verb!bPP-yLB-check.gp!
in the \verb!bPP! directory of \url{https://github.com/PV-314/hypgeom}.

For part~(a), when $b=p^{m}$ is an odd prime power, we have some examples showing
that part~(a) is close to best possible:\\
--$(a,b,t,u)=(685, 71647, 18386, 25)$ and
$y_{1}=18068^{2}=b^{2}/15.724\ldots$,\\
--$(a,b,t,u)=(353, 19207, 4930, 13)$ and
$y_{1}=4844^{2}=b^{2}/15.722\ldots$,\\
--$(a,b,t,u)=(103489,361819267, 91017536, 865)$ and
$y_{1}=90641873^{2}=b^{2}/15.933\ldots$.

From the search that led to the above examples, it seems likely that there are
infinite families with $y_{\pm 1}/b^{2} \rightarrow 1/16$. However, the best
family we found has $y_{-1}/b^{2}$ approaching $1/9$ from below.
For any integer $u \geq 3$, let $a=u-2$, $t=2u^{2}+2$,
$b=3u^{2}+4u+4$, so $d=4u^{2}+8$ and $y_{-1}=\left( u^{2}+u+1 \right)^{2}$.
Here $9y_{-1}=9u^{4}+18u^{3}+27u^{2}+18u+9$, while $b^{2}=9u^{4}+24u^{3}+40u^{2}+32u+16$,
so $9y_{-1}<b^{2}$. Since the polynomial $3u^{2}+4u+4$ is irreducible, there are
conjecturally infinitely many primes of this form.

When $b=2^{m}$, we have only a handful of examples with $\left| y_{\pm 1} \right| <b^{2}/16$
and a perfect square. The worst example is $a=2$, $b=2^{6}$ with $t=8$ and $u=1$.
Here we have $y_{-1}=25=b^{2}/163.84\ldots$

\vspace*{1.0mm}

For part~(b), suppose that $b=2^{\ell}p^{m}$ where $\ell$ and $m$ are positive
integers and $p$ is an odd prime.

For any integer $u>2^{2-\ell/2}$, put $a=4^{\ell-3}u^3-2^{\ell-2}u$, 
$b=8^{\ell-2}u^{4}+4^{\ell-2}u^{2}+2^{\ell}$
and $t=4^{\ell-3}u^{4}+2^{\ell-2}u^{2}+2$. Then
$y_{1}= \left( 2^{2\ell-6}u^{4}+2^{\ell-3}u^{2}+1 \right)^{2}$.
So $y_{1}/b^{2}=2^{4\ell-12}u^{8}/ \left( 2^{6\ell-12} u^{8} \right)+\text{lower-order terms}
=2^{-2\ell}+\cdots$, as required.
Subject to well-known conjectures in analytic number theory, there are infinitely
many $u$ such that $b=2^{\ell}p^{m}$.

For $\ell=2$, we have the family of examples, $(a,b,t,u)=(t-4,8t-12,t, 1)$.
Here $y_{1}=(t-1)^{2}$, so $y_{1}/b^{2}=1/64+\text{lower-order terms}$.

\section{Acknowledgements}

I gratefully thank the referee. Their very careful reading of the manuscript and
their very helpful suggestions certainly improved it in significant and welcome
ways. I also appreciate their
good suggestion to make the code used for this work available via GitHub.

Once again (as with \cite{V5}), I thank Mihai Cipu for the time and effort he
generously spent reading and commenting on an earlier version of this work, as
well as for his stimulating discussions.

\section{Data Availability}

No data beyond that present in the paper itself is used in this work.

\section{Conflict of Interests}

None.

\end{document}